# MAXIMUM SMOOTHED LIKELIHOOD ESTIMATION AND SMOOTHED MAXIMUM LIKELIHOOD ESTIMATION IN THE CURRENT STATUS MODEL

By Piet Groeneboom, Geurt Jongbloed and Birgit I. Witte

*Delft University of Technology, Delft University of Technology and Delft University of Technology*

We consider the problem of estimating the distribution function, the density and the hazard rate of the (unobservable) event time in the current status model. A well studied and natural nonparametric estimator for the distribution function in this model is the nonparametric maximum likelihood estimator (MLE). We study two alternative methods for the estimation of the distribution function, assuming some smoothness of the event time distribution. The first estimator is based on a maximum smoothed likelihood approach. The second method is based on smoothing the (discrete) MLE of the distribution function. These estimators can be used to estimate the density and hazard rate of the event time distribution based on the plug-in principle.

## 1. Introduction.

In survival analysis, one is interested in the distribution of the time it takes before a certain event (failure, onset of a disease) takes place. Depending on exactly what information is obtained on the time $X$ and the precise assumptions imposed on its distribution function $F_0$, many estimators for $F_0$ have been defined and studied in the literature.

When a sample of $X_i$'s is directly and completely observed, one can estimate $F_0$ under various assumptions. In the parametric approach, one assumes $F_0$ to belong to a parametric class of distributions, e.g., the exponential- or Weibull distributions. Then estimating $F_0$ boils down to estimating a finite-dimensional parameter and a variety of classical point estimation procedures can be used to do this. If one wishes to estimate $F_0$ fully nonparametrically, so without assuming any properties of $F_0$ other than the basic properties of distribution functions, the empirical distribution function $\mathbb{F}_n$ of $X_1, \dots, X_n$ is a natural candidate to use. If the distribution function is known to have a continuous derivative $f_0$ w.r.t. Lebesgue measure, one can use kernel estimators [see, e.g., Silverman (1986)] or wavelet methods [see, e.g., Donoho and Johnstone (1995)] for estimating $f_0$. Finally, in case $F_0$ is known to satisfy a certain shape constraint as concavity or convex-concavity on $[0, \infty)$, a shape-constrained estimator for $F_0$ can be used. Problems of this type were considered in, e.g., Bickel and Fan (1996), Groeneboom, Jongbloed and Wellner (2002) and Dümbgen and Rufibach (2009).

However, in many cases the variable $X$ is not observed completely, due to some sort of censoring. Parametric inference in such situations is often not really different from that based on exactly observed $X_i$'s. The parametric model for $X$ basically transforms to a parametric model for the observable data and the usual methods for parametric point estimation can be used to estimate $F_0$. For various types of censoring, also nonparametric estimators have been proposed. In the context of right-censoring, the Kaplan–Meier estimator [see Kaplan and Meier (1958)] is the (nonparametric) maximum likelihood estimator of $F_0$. It maximizes the likelihood of the observed data over all distribution functions, without any additional constraints. Density estimators also exist in this setting, see, e.g., Marron and Padgett (1987). Huang and Zhang (1994) consider the MLE for estimating $F_0$ and its density in this setting under the assumption that $F_0$ is concave on $[0, \infty)$.

The type of censoring we focus on in this paper, is interval censoring, case I. The model for this type of observations is also known as the *current status model*. In this model, a censoring variable $T$, independent of $X$, is observed as well as a variable $\Delta = 1_{\{X \le T\}}$, indicating whether the (unobservable) $X$ lies to the left or to the right of the observed $T$. For this model, the (nonparametric) maximum likelihood estimator is studied in Groeneboom and Wellner (1992). This estimator is discrete and is therefore not suitable for estimating the density $f_0$, the hazard rate $\lambda_0 = f_0/(1 - F_0)$ or the transmission potential which depend on the hazard rate $\lambda_0$ studied in Keiding (1991). An estimator that can be used to estimate these quantities is the maximum likelihood estimator studied by Dümbgen, Freitag-Wolf and Jongbloed (2006) under the constraint that $F$ is concave or convex-concave.

In this paper, we study two likelihood based estimators for $F_0$ (and its density $f_0$ and hazard rate $\lambda_0$) based on interval censored data from $F_0$ under the assumption that $F_0$ is continuously differentiable. The first estimator we study is a so-called maximum smoothed likelihood estimator (MSLE) as studied by Eggermont and LaRiccia (2001) in the context of monotone and unimodal density estimation. It is a general likelihood-based $M$-estimator that will turn out to be smooth automatically. The second estimator we consider, the smoothed maximum likelihood estimator (SMLE), is obtained by convolving the (discrete) MLE of Groeneboom and Wellner (1992) with a smoothing kernel. These different methods result in different but related estimators. Analyzing the pointwise asymptotics shows that only the biases of these estimators differ while the variances are equal. We cannot say that one estimator is uniformly superior to the other. In a somewhat analogous way, Mammen (1991) studies the differences between the efficiencies of smoothing of isotonic estimates and isotonizing smooth estimates. This also does not produce a clear "winner."

The outline of this paper is as follows. In Section 2, we introduce the current status model and review some results needed in the sequel. The MSLE $\hat{F}_n^{\mathrm{MS}}$ for $F_0$ based on current status data is introduced and characterized in Section 3. Moreover, asymptotic results are derived for $\hat{F}_n^{\mathrm{MS}}$ as well as its density $\hat{f}_n^{\mathrm{MS}}$ and hazard rate $\hat{\lambda}_n^{\mathrm{MS}}$, showing that the rate of convergence of $\hat{F}_n^{\mathrm{MS}}$ is faster than the rate of convergence of the MLE. In Section 4, the SMLE for $F_0$, $f_0$ and $\lambda_0$ are introduced and their asymptotic properties derived. The resulting asymptotic distributions are very similar to the asymptotic distributions of the MSLE. In Section 5, we briefly address the problem of bandwidth selection in practice. We also apply these methods to a data set on hepatitis A from Keiding (1991). Technical proofs and lemmas can be found in the Appendix.

## 2. The current status model.

Consider an i.i.d. sequence $X_1, X_2, \ldots$ with distribution $F_0$ on $[0,\infty)$ and independent of this an i.i.d. sequence $T_1, T_2, \ldots$ from a distribution $G$ with Lebesgue density $g$ on $[0,\infty)$. Based on these sequences, define $Z_i = (T_i, 1_{\{X_i \le T_i\}}) =: (T_i, \Delta_i)$. Then $Z_1, Z_2, \ldots$ are i.i.d. and have density $f_Z$ with respect to the product of Lebesgue- and counting measure on $[0,\infty)\times\{0,1\}$:

$$
\begin{aligned}
f_Z(t,\delta) &= g(t)\{\delta F_0(t) + (1-\delta)(1-F_0(t))\} \\
&= \delta g_1(t) + (1-\delta) g_0(t).
\end{aligned}
\tag{2.1}
$$

One usually says that the $X_i$'s take their values in the hidden space $[0,\infty)$ and the $Z_i$ take their values in the observation space $[0,\infty)\times\{0,1\}$.

Let $\mathbb{P}_n$ be the empirical distribution of $Z_1, \ldots, Z_n$. Writing down the log likelihood as a function of $F$ and dividing by $n$, we get

$$
l(F) = \int \{\delta \log F(t) + (1-\delta)\log(1-F(t))\}\, d\mathbb{P}_n(t,\delta).
\tag{2.2}
$$

Here, we ignore a term in the log likelihood that does not depend on the distribution function $F$.

In Groeneboom and Wellner (1992), it is shown that the (nonparametric) maximum likelihood estimator (MLE) is well defined as maximizer of (2.2) over all distribution functions and that it can be characterized as the left derivative of the greatest convex minorant of a cumulative sum diagram. To be precise, the observed time points $T_i$ are ordered in increasing order, yielding $T_{(1)} < T_{(2)} < \cdots < T_{(n)}$, and the $\Delta$ associated with $T_{(i)}$ is denoted by $\Delta_{(i)}$. Then the cumulative sum diagram consisting of the points

$$P_0 = (0,0), \qquad P_i = \left(\frac{i}{n}, \frac{1}{n}\sum_{j=1}^{i} \Delta_{(j)}\right)$$

is constructed. Having determined the greatest convex minorant of this diagram, $\hat{F}_n(T_{(i)})$ is given by the left derivative of this minorant, evaluated at the point $P_i$. At other points it is defined by right continuity. Denoting by $\mathbb{G}_n$ the empirical distribution function of the $T_i$'s and by $\mathbb{G}_{n,1}$ the empirical subdistribution function of the $T_i$'s with $\Delta_i = 1$, observe that for $0 \le i \le n$, $P_i = (\mathbb{G}_n(T_{(i)}), \mathbb{G}_{n,1}(T_{(i)}))$. Also note that $\hat{F}_n$ is a step function of which the set of jump points $\{\tau_1, \ldots, \tau_m\}$ is a subset of the set $\{T_i : 1 \le i \le n\}$.

Groeneboom and Wellner (1992) show that this MLE is a consistent estimator of $F_0$, and prove that under some local smoothness assumptions, for $t > 0$ fixed, $n^{1/3}(\hat{F}_n(t) - F_0(t))$ has the so-called Chernoff distribution as limiting distribution. If $F_0$ and $G$ are assumed to satisfy conditions (F.1) and (G.1) below Groeneboom and Wellner (1992) also prove (see their Lemma 5.9 and page 120)

$$\|F_0 - \hat{F}_n\|_\infty = \mathcal{O}_p(n^{-1/3}\log n), \tag{2.3}$$

$$\max_{1 \le i \le m} |\tau_{i+1} - \tau_i| = \mathcal{O}_p(n^{-1/3}\log n). \tag{2.4}$$

(F.1) $F_0$ has bounded support $\mathcal{S}_0 = [0, M_0]$ and is strictly increasing on $\mathcal{S}_0$ with density $f_0$, strictly staying away from zero.

(G.1) $G$ has support $\mathcal{S}_G = [0, \infty)$, is strictly increasing on $\mathcal{S}_0$ with density $g$ staying away from zero and $g'$ is bounded on $\mathcal{S}_0$.

From this, it follows that for fixed $t > 0$, any $\nu > 0$ and $\mathcal{I}_t = [t - \nu, t + \nu]$

$$\sup_{u \in \mathcal{I}_t} |F_0(u) - \hat{F}_n(u)| = \mathcal{O}_p(n^{-1/3}\log n), \tag{2.5}$$

$$\max_{i : \tau_i \in \mathcal{I}_t} |\tau_{i+1} - \tau_i| = \mathcal{O}_p(n^{-1/3}\log n). \tag{2.6}$$

If one is willing to assume smoothness on $F_0$ and use this in the estimation procedure, this cube-root-$n$ rate of convergence of the estimator can be improved. The two estimators of $F_0$ we define, do indeed converge at the faster rate $n^{2/5}$.

**3. Maximum smoothed likelihood estimation.** In this section, we define the maximum smoothed likelihood estimator (MSLE) $\hat{F}_n^{\mathrm{MS}}$ for the unknown distribution function $F_0$ of the variable of interest $X$. We characterize this estimator as the derivative of the convex minorant of a function on $\mathbb{R}$ and derive its pointwise asymptotic distribution. Based on $\hat{F}_n^{\mathrm{MS}}$, estimators for the density $f_0$ as well as for the hazard rate $\lambda_0 = f_0/(1-F_0)$ are defined and studied asymptotically.

We start with defining the estimators. Define the empirical subdistribution functions based on the $T_j$'s with $\Delta_j = 0$ and 1, respectively, by

$$\mathbb{G}_{n,i}(t) = \frac{1}{n}\sum_{j=1}^{n} 1_{[0,t]\times\{i\}}(T_j,\Delta_j) \qquad \text{for } i = 0, 1,$$

and note that the empirical distribution of the data $\{Z_j = (T_j,\Delta_j) : 1 \le j \le n\}$ can be expressed as $d\mathbb{P}_n(t,\delta) = \delta\, d\mathbb{G}_{n,1}(t) + (1-\delta)\, d\mathbb{G}_{n,0}(t)$. Let $\hat{G}_{n,1}$ and $\hat{G}_{n,0}$ be smoothed versions of $\mathbb{G}_{n,1}$ and $\mathbb{G}_{n,0}$, respectively (e.g., via kernel smoothing), let $\hat{g}_{n,1}$ and $\hat{g}_{n,0}$ be their densities w.r.t. Lebesgue measure on $[0,\infty)$ and define $d\hat{P}_n(t,\delta) = \delta\, d\hat{G}_{n,1}(t) + (1-\delta)\, d\hat{G}_{n,0}(t)$. This is a smoothed version of the empirical measure $\mathbb{P}_n$, where smoothing is only performed "in the $t$-direction." Following the general approach of Eggermont and LaRiccia (2001), we replace the empirical distribution $\mathbb{P}_n$ in the definition of the log likelihood (2.2) by this smoothed version $\hat{P}_n$, and define the smoothed log likelihood on the class of all distribution functions by

$$\begin{aligned} l^S(F) &= \int \{\delta \log F(t) + (1-\delta)\log(1-F(t))\}\, d\hat{P}_n(t,\delta) \\ &= \int \log(1-F(t))\, d\hat{G}_{n,0}(t) + \int \log F(t)\, d\hat{G}_{n,1}(t). \end{aligned} \tag{3.1}$$

The maximizer of the smoothed log likelihood is characterized similarly as the maximizer of the log likelihood. The next theorem makes this precise.

THEOREM 3.1. *Define $\hat{G}_n(t) = \hat{G}_{n,0}(t) + \hat{G}_{n,1}(t)$ for $t \ge 0$ and consider the following parameterized curve in $\mathbb{R}_+^2$, a continuous cumulative sum diagram (CCSD):*

$$t \mapsto (\hat{G}_n(t), \hat{G}_{n,1}(t)), \tag{3.2}$$

*for $t \in [0,\tau]$, with $\tau = \sup\{t \ge 0 : \hat{g}_{n,0}(t) + \hat{g}_{n,1}(t) > 0\}$. Let $\hat{F}_n^{\mathrm{MS}}(t)$ be the right-continuous slope of the lower convex hull of the CCSD (3.2), evaluated at the point with $x$-coordinate $\hat{G}_n(t)$. Then $\hat{F}_n^{\mathrm{MS}}$ is the unique maximizer of (3.1) over the class of all sub-distribution functions. We call $\hat{F}_n^{\mathrm{MS}}$ the maximum smoothed likelihood estimator of $F_0$.*

In the proof of Theorem 3.1, we use the following lemma, a proof of which can be found in the Appendix.

LEMMA 3.2. *Let* $\hat{F}_n^{\mathrm{MS}}$ *be defined as in Theorem 3.1. Then for any distribution function* $F$,

$$\int \log F(t)\, d\hat{G}_{n,1}(t) \le \int \hat{F}_n^{\mathrm{MS}}(t) \log F(t)\, d\hat{G}_n(t)$$

*and*

$$\int \log(1-F(t))\, d\hat{G}_{n,0}(t) \le \int (1-\hat{F}_n^{\mathrm{MS}}(t)) \log(1-F(t))\, d\hat{G}_n(t)$$

*with equality in case* $F = \hat{F}_n^{\mathrm{MS}}$.

PROOF OF THEOREM 3.1. Use the equality part of Lemma 3.2 to rewrite (3.1) as

$$l^S(\hat{F}_n^{\mathrm{MS}}) = \int (\hat{F}_n^{\mathrm{MS}}(t) \log \hat{F}_n^{\mathrm{MS}}(t) + (1-\hat{F}_n^{\mathrm{MS}}(t)) \log(1-\hat{F}_n^{\mathrm{MS}}(t)))\, d\hat{G}_n(t).$$

By the inequality part of Lemma 3.2, we get for each distribution function $F$ that

$$l^S(F) \le \int \hat{F}_n^{\mathrm{MS}}(t) \log F(t)\, d\hat{G}_n(t) + \int (1-\hat{F}_n^{\mathrm{MS}}(t)) \log(1-F(t))\, d\hat{G}_n(t).$$

Now note, using the convention $0 \cdot \infty = 0$, that for all $p, p' \in [0,1]$

$$p \log p' + (1-p)\log(1-p') \le p \log p + (1-p) \log(1-p). \tag{3.3}$$

This implies that $l^S(F) \le l^S(\hat{F}_n^{\mathrm{MS}})$, i.e., $l^S$ is maximal for $\hat{F}_n^{\mathrm{MS}}$.

For uniqueness, note that inequality (3.3) is strict whenever $p' \ne p$. The last step in the preceding argument then shows that $l^S(F) < l^S(\hat{F}_n^{\mathrm{MS}})$, unless $F = \hat{F}_n^{\mathrm{MS}}$ a.e. w.r.t. the measure $d\hat{G}_n$. It could be that $d\hat{G}_n$ has no mass on $[a,b]$ for some $a<b$, i.e., $(\hat{G}_n(t), \hat{G}_{n,1}(t)) = (\hat{G}_n(a), \hat{G}_{n,1}(a))$ for all $t \in [a,b]$. This means that $\hat{F}_n^{\mathrm{MS}}$ is constant on $[a,b]$. Furthermore, it holds that $F(a) = \hat{F}_n^{\mathrm{MS}}(a)$ and $F(b) = \hat{F}_n^{\mathrm{MS}}(b)$, implying that $F$ is also constant and equal to $\hat{F}_n^{\mathrm{MS}}$ on $[a,b]$ a.e. w.r.t. the Lebesgue measure on $[0,\infty)$. Hence, $l^S(F) < l^S(\hat{F}_n^{\mathrm{MS}})$ unless $F = \hat{F}_n^{\mathrm{MS}}$. □

We assume the estimators $\hat{G}_{n,i}$ are continuously differentiable, hence, $\hat{F}_n^{\mathrm{MS}}$ is continuous and its derivative exists. So we can define the maximum smoothed likelihood estimators for $f_0$ and $\lambda_0$ by

$$\hat{f}_n^{\mathrm{MS}}(t) = \frac{d}{du}\hat{F}_n^{\mathrm{MS}}(u)\Big|_{u=t}, \qquad \hat{\lambda}_n^{\mathrm{MS}}(t) = \frac{\hat{f}_n^{\mathrm{MS}}(t)}{1-\hat{F}_n^{\mathrm{MS}}(t)} \tag{3.4}$$

for $t > 0$ such that $\hat{F}_n^{\mathrm{MS}}(t) < 1$.

In Theorem 3.1 no particular choice for $\hat{G}_{n,0}$ and $\hat{G}_{n,1}$ was made. For what follows, we define these estimators explicitly as kernel smoothed versions of $\mathbb{G}_{n,0}$ and $\mathbb{G}_{n,1}$. Let $k$ be a probability density satisfying condition (K.1).

(K.1) The probability density $k$ has support $[-1, 1]$, is symmetric and twice continuously differentiable on $\mathbb{R}$.

Note that condition (K.1) implies that $m_2(k) = \int u^2 k(u)\, du < \infty$.

Let $K$ be the distribution function with density $k$, i.e., $K(t) = \int_{-\infty}^{t} k(u)\, du$, $k'$ be the derivative of $k$ and $h > 0$ be a smoothing parameter (depending on $n$). Then we use the following notation for the scaled version of $K$, $k$ and $k'$

$$K_h(u) = K(u/h), \qquad k_h(u) = \frac{1}{h} k(u/h) \quad \text{and} \quad k_h'(u) = \frac{1}{h^2} k'(u/h). \tag{3.5}$$

For $i = 0, 1$ let

$$\hat{g}_{n,i}(t) = \int k_h(t-u)\, d\mathbb{G}_{n,i}(u)$$

be kernel (sub-density) estimates based on the observations $T_j$ for which $\Delta_j = i$, and let $\hat{g}_n(t) = \hat{g}_{n,1}(t) + \hat{g}_{n,0}(t)$. Also define the associated (sub-) distribution functions

$$\hat{G}_{n,i}(t) = \int_{[0,t]} \hat{g}_{n,i}(u)\, du, \qquad \text{for } i = 0, 1, \quad \text{and} \quad \hat{G}_n(t) = \int_{[0,t]} \hat{g}_n(u)\, du.$$

Because $X \geq 0$, we can expect inconsistency problems for the kernel density and density derivative estimators at zero. In order to prevent those, we modify the definition of $\hat{g}_{n,i}$ for $t < h$. To be precise, we define

$$\hat{g}_{n,i}(t) = \int \frac{1}{h} k^{\beta}\left(\frac{t-u}{h}\right) d\mathbb{G}_{n,i}(u), \qquad 0 \leq t \leq h,$$

for $\beta = t/h$ where the so-called boundary kernel $k^{\beta}$ is defined by

$$k^{\beta}(u) = \frac{\nu_{2,\beta}(k) - \nu_{1,\beta}(k) u}{\nu_{0,\beta}(k)\nu_{2,\beta}(k) - \nu_{1,\beta}(k)^2} k(u) 1_{(-1,\beta)}(u)$$

$$\text{with } \nu_{i,\beta}(k) = \int_{-1}^{\beta} u^i k(u)\, du, i = 0, 1, 2.$$

Let the estimators $\hat{g}_{n,i}'$ be the derivatives of $\hat{g}_{n,i}$, for $i = 0, 1$. There are other ways to correct the kernel estimator near the boundary, see, e.g., Schuster (1985) or Jones (1993). However, simulations show that the results are not much influenced by the used boundary correction method.

Having made these choices for the smoothed empirical distribution $\hat{P}_n$, let us return to the MSLE. It is the maximizer of $l^S$ over the class of all distribution functions. One could also maximize $l^S$ over the bigger class of all functions, maximizing the integrand of (3.1) for each $t$ separately. This results in

$$\hat F_n^{\mathrm{naive}}(t) = \frac{\hat g_{n,1}(t)}{\hat g_n(t)}, \qquad \hat f_n^{\mathrm{naive}}(t) = \frac{\hat g_n(t)\hat g'_{n,1}(t) - \hat g'_n(t)\hat g_{n,1}(t)}{\hat g_n(t)^2}, \tag{3.6}$$

where

$$\hat g'_n(t) = \hat g'_{n,0}(t) + \hat g'_{n,1}(t). \tag{3.7}$$

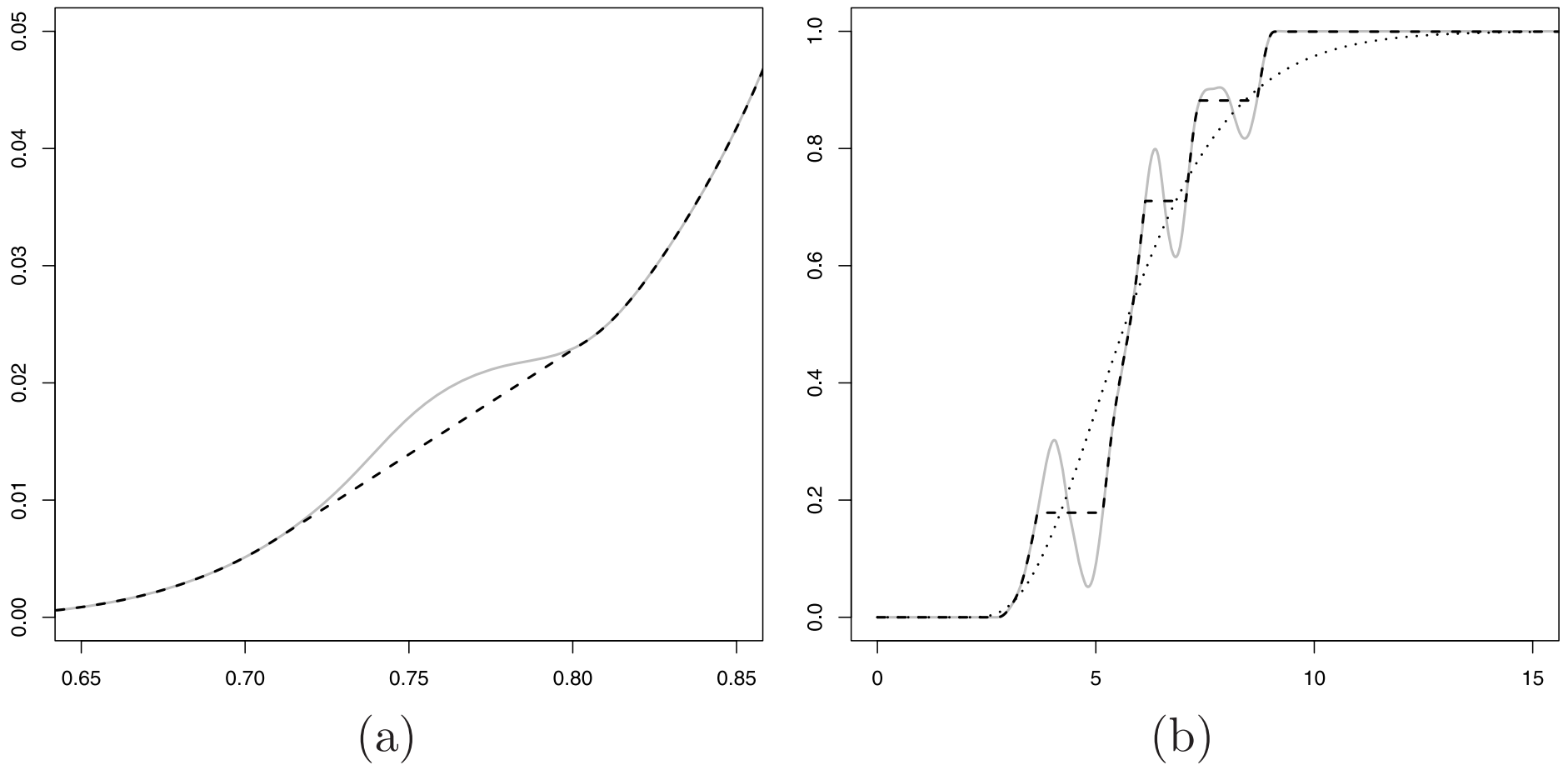

FIG. 1. *A part of the CCSD, its lower convex hull and the estimates* $\hat F_n^{\mathrm{naive}}$ *and* $\hat F_n^{\mathrm{MS}}$ *for* $F_0$ *based on simulated data, with* $n = 500$. (a) *Part of the CCSD (grey line) and its lower convex hull (dashed line);* (b) *estimates* $\hat F_n^{\mathrm{naive}}$ *(grey line) and* $\hat F_n^{\mathrm{MS}}$ *(dashed line) of* $F_0$ *(dotted line).*

We call these naive estimators, since $\hat f_n^{\mathrm{naive}}$ might take negative values, meaning that $\hat F_n^{\mathrm{naive}}$ decreases locally.

Figure 1(a) shows a part of the CCSD defined in (3.2) and its lower convex hull. Figure 1(b) shows the naive estimator $\hat F_n^{\mathrm{naive}}$ (the grey line), the MSLE $\hat F_n^{\mathrm{MS}}$ and the true distribution for a simulation of size 500. The unknown distribution of the variable $X$ is taken to be a shifted Gamma(4) distribution, i.e., $f_0(x) = \frac{(x-2)^3}{3!}\exp(-(x-2))1_{[2,\infty)}(x)$, and the censoring variable $T$ has an exponential distribution with mean 3, i.e., $g(t) = \frac{1}{3}\exp(-t/3)1_{[0,\infty)}$. For the kernel density, we took the triweight kernel $k(t) = \frac{35}{32}(1-t^2)^3 1_{[-1,1]}(t)$ and as bandwidth $h = 0.7$. This picture shows that the estimator $\hat F_n^{\mathrm{MS}}$ is the isotonic version of the estimator $\hat F_n^{\mathrm{naive}}$.

The next theorem shows that for appropriately chosen $h$, the naive estimator $\hat F_n^{\mathrm{naive}}$ will be monotonically increasing on big intervals with probability converging to one as $n$ tends to infinity if $F_0$ and $G$ satisfy conditions (F.1) and (G.1).

THEOREM 3.3. *Assume $F_0$ and $G$ satisfy conditions* (F.1) *and* (G.1). *Let $\hat{g}_n$ and $\hat{g}_{n,1}$ be kernel estimators for $g$ and $g_1$ with kernel density $k$ satisfying condition* (K.1). *Let $h = cn^{-\alpha}$ $(c > 0)$ be the bandwidth used in the definition of $\hat{g}_n$ and $\hat{g}_{n,1}$. Then for all $0 < m < M < M_0$ and $\alpha \in (0, 1/3)$ the following holds*

$$P(\hat{F}_n^{\text{naive}} \textit{ is monotonically increasing on } [m, M]) \longrightarrow 1. \tag{3.8}$$

Note that this theorem as it stands does not imply that $\hat{F}_n^{\text{MS}}(t) = \hat{F}_n^{\text{naive}}(t)$ on $[m, M]$ with probability tending to one. Some additional control on the behavior of $\hat{F}_n^{\text{naive}}$ on $[0, m)$ and $(M, M_0]$ is needed. The proof of the corollary below makes this precise.

COROLLARY 3.4. *Under the assumptions of Theorem 3.3, it holds that for all $0 < m < M < M_0$ and $\alpha \in (0, 1/3)$,*

$$P(\hat{F}_n^{\text{naive}}(t) = \hat{F}_n^{\text{MS}}(t) \textit{ for all } t \in [m, M]) \longrightarrow 1. \tag{3.9}$$

*Consequently, for all $t > 0$ the asymptotic distributions of $\hat{F}_n^{\text{MS}}(t)$ and $\hat{F}_n^{\text{naive}}(t)$ are the same.*

In van der Vaart and van der Laan (2003), a result similar to our Corollary 3.4 is proved for smooth monotone density estimators. The kernel estimator is compared with an isotonized version of this estimator. Their proof is based on a so-called switch-relation relating the derivative of the convex minorant of a function to that of an argmax function. The direct argument we use to prove Corollary 3.4 furnishes an alternative way to prove their result.

By Corollary 3.4, the estimators $\hat{F}_n^{\text{MS}}(t)$ and $\hat{F}_n^{\text{naive}}(t)$ have the same asymptotic distribution. The same holds for $\hat{f}_n^{\text{MS}}(t)$ and $\hat{f}_n^{\text{naive}}(t)$ as well as for $\hat{\lambda}_n^{\text{MS}}(t)$ and $\hat{\lambda}_n^{\text{naive}}(t)$. The pointwise asymptotic distribution of $\hat{F}_n^{\text{naive}}(t)$ follows easily from the Lindeberg–Feller central limit theorem and the delta method. The resulting pointwise asymptotic normality of both $\hat{F}_n^{\text{MS}}(t)$ and $\hat{F}_n^{\text{naive}}(t)$ is stated in the next theorem.

THEOREM 3.5. *Assume $F_0$ and $G$ satisfy conditions* (F.1) *and* (G.1). *Fix $t > 0$ such that $f_0''$ and $g''$ exist and are continuous at $t$ and $g(t)f_0'(t) + 2f_0(t)g'(t) \neq 0$. Let $h = cn^{-1/5}$ $(c > 0)$ be the bandwidth used in the definition of $\hat{g}_n$ and $\hat{g}_{n,1}$. Then*

$$n^{2/5}(\hat{F}_n^{\text{MS}}(t) - F_0(t)) \rightsquigarrow \mathcal{N}(\mu_{F,\text{MS}}, \sigma^2_{F,\text{MS}}),$$

*where*

$$\mu_{F,\mathrm{MS}} = \frac{1}{2}c^2 m_2(k)\left\{f_0'(t) + 2\frac{f_0(t)g'(t)}{g(t)}\right\},$$

$$\sigma^2_{F,\mathrm{MS}} = c^{-1}\frac{F_0(t)(1-F_0(t))}{g(t)}\int k(u)^2\,du.$$

*This also holds if we replace* $\hat{F}_n^{\mathrm{MS}}$ *by* $\hat{F}_n^{\mathrm{naive}}$.

*For fixed* $t>0$, *the asymptotically MSE-optimal bandwidth* $h$ *for* $\hat{F}_n^{\mathrm{MS}}(t)$ *is given by* $h_{n,F,\mathrm{MS}} = c_{F,\mathrm{MS}} n^{-1/5}$, *where*

$$\begin{aligned} c_{F,\mathrm{MS}} &= \left\{\frac{F_0(t)(1-F_0(t))}{g(t)}\int k(u)^2\,du\right\}^{1/5} \\ &\quad\times\left\{m_2^2(k)\left\{f_0'(t) + 2\frac{f_0(t)g'(t)}{g(t)}\right\}^2\right\}^{-1/5}. \end{aligned} \tag{3.10}$$

PROOF. For fixed $c>0$, the asymptotic distribution of $\hat{F}_n^{\mathrm{naive}}$ follows immediately by applying the delta method with $\varphi(u,v) = v/u$ to the first result in Lemma A.3. By Corollary 3.4, this also gives the asymptotic distribution of $\hat{F}_n^{\mathrm{MS}}$.

To obtain the bandwidth which minimizes the asymptotic mean squared error (aMSE) we minimize

$$\begin{aligned} \mathrm{aMSE}(\hat{F}_n^{\mathrm{MS}}, c) &= \frac{1}{4}c^4 m_2^2(k)\left\{f_0'(t) + 2\frac{f_0(t)g'(t)}{g(t)}\right\}^2 \\ &\quad + c^{-1}\frac{F_0(t)(1-F_0(t))}{g(t)}\int k(u)^2\,du \end{aligned}$$

with respect to $c$. This yields (3.10). □

REMARK 3.1. In case $g(t)f_0'(t) + 2f_0(t)g'(t) = 0$, the optimal rate of $h_{n,F,\mathrm{MS}}$ is $n^{-1/9}$ resulting in a rate of convergence $n^{-4/9}$ for $\hat{F}_n^{\mathrm{MS}}$. This is in line with results for other kernel smoothers in case of vanishing first-order bias terms.

The pointwise asymptotic distributions of $\hat{f}_n^{\mathrm{MS}}(t)$ and $\hat{f}_n^{\mathrm{naive}}(t)$ also follow from the Lindeber–Feller central limit theorem and the delta method.

THEOREM 3.6. *Consider* $\hat{f}_n^{\mathrm{MS}}$ *as defined in (3.4) and assume* $F_0$ *and* $G$ *satisfy conditions* (F.1) *and* (G.1). *Fix* $t>0$ *such that* $f_0^{(3)}$ *and* $g^{(3)}$ *exist and are continuous at* $t$. *Let* $h = cn^{-1/7}$ ($c>0$) *be the bandwidth used to define* $\hat{F}_n^{\mathrm{MS}}$. *Then*

$$n^{2/7}(\hat{f}_n^{\mathrm{MS}}(t) - f_0(t)) \rightsquigarrow \mathcal{N}(\mu_{f,\mathrm{MS}}, \sigma^2_{f,\mathrm{MS}}),$$

*where*

$$\begin{aligned}\mu_{f,\mathrm{MS}} &= \frac{1}{2}c^2 m_2(k)\left(f_0''(t) + 2\frac{g''(t)f_0(t)+g'(t)f_0'(t)}{g(t)} - 2\frac{g'(t)^2 f_0(t)}{g(t)^2}\right)\\ &=: \frac{1}{2}c^2 m_2(k)q(t),\\ \sigma^2_{f,\mathrm{MS}} &= \frac{F_0(t)(1-F_0(t))}{c^3 g(t)}\int k'(u)^2\,du\end{aligned}$$

*for $t$ such that $q(t)\neq 0$. This also holds if we replace $\hat f_n^{\mathrm{MS}}$ by $\hat f_n^{\mathrm{naive}}$.*

*For fixed $t>0$, the aMSE-optimal bandwidth $h$ for $\hat f_n^{\mathrm{MS}}(t)$ is given by $h_{n,f,\mathrm{MS}} = c_{f,\mathrm{MS}} n^{-1/7}$, where*

$$c_{f,\mathrm{MS}} = \left\{3\frac{F_0(t)(1-F_0(t))}{g(t)}\int k'(u)^2\,du\right\}^{1/7}\{m_2^2(k)q^2(t)\}^{-1/7}. \tag{3.11}$$

Proof. Write $\hat g_n(t) = g(t) + R_n(t)$ and $\hat g_{n,1}(t) = g_1(t) + R_{n,1}(t)$, so

$$\begin{aligned}&n^{2/7}(\hat f_n^{\mathrm{naive}}(t) - f_0(t))\\ &\quad = n^{2/7}\left(\frac{g(t)\hat g_{n,1}'(t) - \hat g_n'(t)g_1(t)}{g(t)^2} - \frac{g(t)g_1'(t) - g'(t)g_1(t)}{g(t)^2}\right) + T_n(t)\end{aligned}$$

for

$$\begin{aligned}T_n(t) &= n^{2/7}\frac{[g(t)+R_n(t)]\hat g_{n,1}'(t) - \hat g_n'(t)[g_1(t)+R_{n,1}(t)]}{[g(t)+R_n(t)]^2}\\ &\quad - n^{2/7}\frac{g(t)\hat g_{n,1}'(t) - \hat g_n'(t)g_1(t)}{g(t)^2}\\ &= n^{2/7}\frac{R_n(t)\hat g_{n,1}'(t) - \hat g_n'(t)R_{n,1}(t)}{[g(t)+R_n(t)]^2}\\ &\quad - n^{2/7}(g(t)\hat g_{n,1}'(t) - \hat g_n'(t)g_1(t))\frac{R_n(t)(2g(t)+R_n(t))}{g(t)^2[g(t)+R_n(t)]^2}.\end{aligned}$$

Applying the delta method with $\varphi(u,v) = (g(t)v - g_1(t)u)/g(t)^2$ to the last result in Lemma A.3 gives that

$$n^{2/7}\left(\frac{g(t)\hat g_{n,1}'(t) - \hat g_n'(t)g_1(t)}{g(t)^2} - \frac{g(t)g_1'(t) - g'(t)g_1(t)}{g(t)^2}\right) \rightsquigarrow \mathcal{N}(\mu_1, \sigma^2_{f,\mathrm{MS}})$$

for

$$\mu_1 = \frac{1}{2}c^2 m_2(k)\left(f_0''(t) + 3\frac{g''(t)f_0(t)+g'(t)f_0'(t)}{g(t)}\right). \tag{3.12}$$

By Lemma A.3 $n^{2/7}R_n(t) \xrightarrow{\mathcal{P}} \frac{1}{2}c^2 m_2(k)g''(t)$ and $n^{2/7}R_{n,1}(t) \xrightarrow{\mathcal{P}} \frac{1}{2}c^2 \times m_2(k)g_1''(t)$, so by the consistency of $\hat{g}_n'$ and $\hat{g}_{n,1}'$, see Lemma A.2, and the continuous mapping theorem we have

$$\begin{aligned} T_n(t) &\xrightarrow{\mathcal{P}} \frac{1}{2}c^2 m_2(k)\frac{g''(t)g_1'(t) - g'(t)g_1''(t)}{g(t)^2} \\ &\qquad - \frac{1}{2}c^2 m_2(k)(g(t)g_1'(t) - g'(t)g_1(t))\frac{2g''(t)g(t)}{g(t)^4} \\ &= -\frac{1}{2}c^2 m_2(k)\left(2\frac{g'(t)^2 f_0(t)}{g(t)^2} + \frac{g''(t)f_0(t) + g'(t)f_0'(t)}{g(t)}\right) = \mu_2. \end{aligned}$$

Hence, we have that

$$n^{2/7}(\hat{f}_n^{\mathrm{naive}}(t) - f_0(t)) \rightsquigarrow \mathcal{N}(\mu_{f,\mathrm{MS}}, \sigma^2_{f,\mathrm{MS}})$$

for $\mu_{f,\mathrm{MS}} = \mu_1 + \mu_2$. By Corollary 3.4, this also gives the asymptotic distribution of $\hat{f}_n^{\mathrm{MS}}$.

The optimal $c$ given in (3.11) is obtained by minimizing

$$\mathrm{aMSE}(\hat{f}_n^{\mathrm{MS}}, c) = \frac{1}{4}c^4 m_2^2(k)q^2(t) + c^{-3}\frac{F_0(t)(1 - F_0(t))}{g(t)}\int k'(u)^2\,du. \qquad \square$$

COROLLARY 3.7. *Consider* $\hat{\lambda}_n^{\mathrm{MS}}$ *of* $\lambda_0$ *as defined in (3.4) and let* $h = cn^{-1/7}$ $(c > 0)$ *be the bandwidth used to compute it. Assume* $F_0$ *and* $G$ *satisfy conditions* (F.1) *and* (G.1). *Fix* $t > 0$ *such that* $F_0(t) < 1$ *and* $f_0^{(3)}$ *and* $g^{(3)}$ *exist and are continuous at* $t$. *Then*

$$n^{2/7}(\hat{\lambda}_n^{\mathrm{MS}}(t) - \lambda_0(t)) \rightsquigarrow \mathcal{N}(\mu_{\lambda,\mathrm{MS}}, \sigma^2_{\lambda,\mathrm{MS}}),$$

*where*

$$\begin{aligned} \mu_{\lambda,\mathrm{MS}} &= \frac{1}{2}c^2 m_2(k)\frac{1}{1 - F_0(t)}\left(f_0''(t) + 2\frac{g''(t)f_0(t) + g'(t)f_0'(t)}{g(t)} - 2\frac{g'(t)^2 f_0(t)}{g(t)^2}\right) \\ &\qquad + \frac{1}{2}c^2 m_2(k)\frac{f_0(t)}{(1 - F_0(t))^2}\left(f_0'(t) + \frac{2g'(t)f_0(t)}{g(t)}\right) = \frac{1}{2}c^2 m_2(k)r(t), \\ \sigma^2_{\lambda,\mathrm{MS}} &= \frac{F_0(t)}{c^3 g(t)(1 - F_0(t))}\int k'(u)^2\,du \end{aligned}$$

*for* $t$ *such that* $r(t) \neq 0$. *This also holds if we replace* $\hat{\lambda}_n^{\mathrm{MS}}$ *by* $\hat{\lambda}_n^{\mathrm{naive}}$.

*For fixed* $t > 0$ *the aMSE-optimal bandwidth* $h$ *for* $\hat{\lambda}_n^{\mathrm{MS}}(t)$ *is given by* $h_{n,\lambda,\mathrm{MS}} = c_{\lambda,\mathrm{MS}} n^{-1/7}$, *where*

$$c_{\lambda,\mathrm{MS}} = \left\{\frac{F_0(t)}{g(t)(1 - F_0(t))}\int k'(u)^2\,du\right\}^{1/7}\{m_2^2(k)r^2(t)\}^{-1/7}. \tag{3.13}$$

PROOF. Write $\hat F_n^{\mathrm{MS}}(t) = F_0(t) + R_n(t)$, then

$$n^{2/7}(\hat\lambda_n^{\mathrm{MS}}(t) - \lambda_0(t)) = \frac{n^{2/7}}{1-F_0(t)}(\hat f_n^{\mathrm{MS}}(t) - f_0(t)) + T_n(t) \tag{3.14}$$

with

$$T_n(t) = n^{2/7}\hat f_n^{\mathrm{MS}}(t)\left(\frac{1}{1-F_0(t)-R_n(t)} - \frac{1}{1-F_0(t)}\right).$$

If $h = cn^{-1/7}$ is the bandwidth for $\hat F_n^{\mathrm{MS}}(t)$, then

$$n^{2/7}\left(\frac{\hat g_{n,1}(t)}{\hat g_n(t)} - \frac{g_1(t)}{g(t)}\right) \xrightarrow{\mathcal{P}} \frac{1}{2}c_2 m_2(k)\left\{f_0'(t) + 2\frac{f_0(t)g'(t)}{g(t)}\right\} = \mu_{F,\mathrm{MS}}$$

by Lemma A.3 and the delta method. This implies that $n^{2/7}R_n(t) \xrightarrow{\mathcal{P}} \mu_{F,\mathrm{MS}}$ and

$$T_n(t) = n^{2/7}\hat f_n^{\mathrm{MS}}(t)\frac{R_n(t)}{(1-F_0(t))(1-F_0(t)-R_n(t))} \xrightarrow{\mathcal{P}} \frac{f_0(t)}{(1-F_0(t))^2}\mu_{F,\mathrm{MS}}.$$

Since we also have that

$$\frac{n^{2/7}}{1-F_0(t)}(\hat f_n^{\mathrm{MS}}(t) - f_0(t)) \rightsquigarrow \mathcal{N}\left(\frac{\mu_{f,\mathrm{MS}}}{1-F_0(t)}, \frac{\sigma^2_{f,\mathrm{MS}}}{(1-F_0(t))^2}\right)$$

we get that $\mu_{\lambda,\mathrm{MS}} = \mu_{f,\mathrm{MS}}/(1-F_0(t)) + \mu_{F,\mathrm{MS}} f_0(t)/(1-F_0(t))^2$.

The optimal $c$ given in (3.13) is obtained by minimizing

$$\mathrm{aMSE}(\hat\lambda_n^{\mathrm{MS}}, c) = \frac{1}{4}c^4 m_2^2(k) r^2(t) + c^{-3}\frac{F_0(t)}{g(t)(1-F_0(t))}\int k'(u)^2\,du. \qquad \square$$

**4. Smoothed maximum likelihood estimation.** In the previous section, we started smoothing the empirical distribution of the observed data, and used that probability measure instead of the empirical distribution function in the definition of the log likelihood. In this section, we consider an estimator that is obtained by smoothing the MLE (see Section 2). Recall the definitions of the scaled versions of $K$, $k$ and $k'$, given in (3.5)

$$K_h(u) = K(u/h), \qquad k_h(u) = \frac{1}{h}k(u/h) \quad \text{and} \quad k_h'(u) = \frac{1}{h^2}k'(u/h).$$

Define the SMLE $\hat F_n^{\mathrm{SM}}$ for $F_0$ by

$$\hat F_n^{\mathrm{SM}}(t) = \int K_h(t-u)\,d\hat F_n(u).$$

Similarly, define the SMLE $\hat f_n^{\mathrm{SM}}$ for $f_0$ and the SMLE $\hat\lambda_n^{\mathrm{SM}}$ of $\lambda_0$ by

$$\hat f_n^{\mathrm{SM}}(t) = \int k_h(t-u)\,d\hat F_n(u) \quad \text{and} \quad \hat\lambda_n^{\mathrm{SM}}(t) = \hat f_n^{\mathrm{SM}}(t)/(1-\hat F_n^{\mathrm{SM}}(t)).$$

In this section, we derive the pointwise asymptotic distributions for these estimators. First, we rewrite the estimators $\hat F_n^{\mathrm{SM}}$ and $\hat f_n^{\mathrm{SM}}$.

LEMMA 4.1. *Fix $t > 0$, such that $g(u) > 0$ in a neighborhood of $t$ and define*

$$\psi_{h,t}(u) = \frac{k_h(t-u)}{g(u)}, \qquad \varphi_{h,t}(u) = \frac{k_h'(t-u)}{g(u)}. \tag{4.1}$$

*Then*

$$\int K_h(t-u)\,d(\hat{F}_n - F_0)(u) = -\int \psi_{h,t}(u)(\delta - \hat{F}_n(u))\,dP_0(u,\delta), \tag{4.2}$$

$$\int k_h(t-u)\,d(\hat{F}_n - F_0)(u) = -\int \varphi_{h,t}(u)(\delta - \hat{F}_n(u))\,dP_0(u,\delta). \tag{4.3}$$

PROOF. To see equality (4.2), we rewrite the left-hand side as follows

$$\begin{aligned}
&\int_0^{t+h} K_h(t-u)\,d(\hat{F}_n - F_0)(u) \\
&\quad = \int_0^{t-h} d(\hat{F}_n - F_0)(u) + \int_{t-h}^{t+h} K_h(t-u)\,d(\hat{F}_n - F_0)(u) \\
&\quad = \hat{F}_n(t-h) - F_0(t-h) + K_h(t-u)(\hat{F}_n(u) - F_0(u))\big|_{u=t-h}^{t+h} \\
&\qquad - \int_{t-h}^{t+h} -(\hat{F}_n(u) - F_0(u))k_h(t-u)\,du \\
&\quad = \int_{t-h}^{t+h} \frac{k_h(t-u)}{g(u)}(\hat{F}_n(u) - F_0(u))\,dG(u) \\
&\quad = -\int \psi_{h,t}(u)(\delta - \hat{F}_n(u))\,dP_0(u,\delta).
\end{aligned}$$

Equation (4.3) follows by a similar argument. □

Hence, in determining the asymptotic distribution of the estimators $\hat{F}_n^{\mathrm{SM}}(t)$ and $\hat{f}_n^{\mathrm{SM}}(t)$, we can consider the integrals at the right-hand side of (4.2) and (4.3). The idea of the proof of the asymptotic result for $\hat{F}_n^{\mathrm{SM}}(t)$, given in the next theorem proven in the Appendix, is as follows. By the characterization of the MLE, given in Lemma A.5, we could add the term $d\mathbb{P}_n$ for free in the right-hand side of (4.2) if $\psi_{h,t}$ were piecewise constant. For most choices of $k$ this function $\psi_{h,t}$ is not piecewise constant. Replacing it by an appropriately chosen piecewise constant function results in an additional $\mathcal{O}_p$-term which does not influence the asymptotic distribution. By some more adding and subtracting, resulting in some more $\mathcal{O}_p$-terms, we get that

$$-n^{2/5}\int \psi_{h,t}(u)(\delta - \hat{F}_n(u))\,dP_0(u,\delta)$$

$$= n^{2/5} \int \psi_{h,t}(u)(\delta - F_0(u))\, d(\mathbb{P}_n - P_0)(u,\delta) + \mathcal{O}_p(1)$$

and the pointwise asymptotic distribution follows from the central limit theorem.

THEOREM 4.2. *Assume $F_0$ and $G$ satisfy conditions* (F.1) *and* (G.1). *Fix $t > 0$ such that $f_0'$ is continuous at $t$ and $f_0'(t) \neq 0$. Let $h = cn^{-\alpha}$ $(c > 0)$ be the bandwidth used in the definition of $\hat{F}_n^{\mathrm{SM}}$. Then for $\alpha = 1/5$*

$$n^{2/5}(\hat{F}_n^{\mathrm{SM}}(t) - F_0(t)) \rightsquigarrow \mathcal{N}(\mu_{F,\mathrm{SM}}, \sigma_{F,\mathrm{SM}}^2),$$

*where*

$$\mu_{F,\mathrm{SM}} = \tfrac{1}{2}c^2 m_2(k) f_0'(t), \qquad \sigma_{F,\mathrm{SM}}^2 = \frac{F_0(t)(1-F_0(t))}{cg(t)} \int k(u)^2\, du. \tag{4.4}$$

*For fixed $t > 0$ the aMSE-optimal bandwidth of $h$ for estimating $\hat{F}_n^{\mathrm{SM}}(t)$ is given by $h_{n,F,\mathrm{SM}} = c_{F,\mathrm{SM}} n^{-1/5}$, where*

$$c_{F,\mathrm{SM}} = \left\{ \frac{F_0(t)(1-F_0(t))}{g(t)} \int k(u)^2\, du \right\}^{1/5} \{m_2^2(k) f_0'(t)^2\}^{-1/5}. \tag{4.5}$$

THEOREM 4.3. *Assume $F_0$ and $G$ satisfy conditions* (F.1) *and* (G.1). *Fix $t > 0$ such that $f_0''$ is continuous at $t$ and $f_0''(t) \neq 0$. Let $h = cn^{-1/7}$ $(c > 0)$ be the bandwidth used in the definition of $\hat{f}_n^{\mathrm{SM}}$. Then*

$$n^{2/7}(\hat{f}_n^{\mathrm{SM}}(t) - f_0(t)) \rightsquigarrow \mathcal{N}(\mu_{f,\mathrm{SM}}, \sigma_{f,\mathrm{SM}}^2),$$

*where*

$$\mu_{f,\mathrm{SM}} = \tfrac{1}{2}c^2 m_2(k) f_0''(t), \qquad \sigma_{f,\mathrm{SM}}^2 = \frac{F_0(t)(1-F_0(t))}{c^3 g(t)} \int k'(u)^2\, du.$$

*For fixed $t > 0$ the aMSE-optimal value of $h$ for estimating $\hat{f}_n^{\mathrm{SM}}(t)$ is given by $h_{n,f,\mathrm{SM}} = c_{f,\mathrm{SM}} n^{-1/7}$, where*

$$c_{f,\mathrm{SM}} = \left\{ 3\frac{F_0(t)(1-F_0(t))}{g(t)} \int k'(u)^2\, du \right\}^{1/7} \{m_2^2(k) f_0''(t)^2\}^{-1/7}. \tag{4.6}$$

The proof of this result is similar to the proof of Theorem 4.2, hence it is omitted.

COROLLARY 4.4. *Assume $F_0$ and $G$ satisfy conditions* (F.1) *and* (G.1). *Fix $t > 0$ such that $F_0(t) < 1$, $f_0''$ is continuous in $t$ and $f_0''(t) \neq 0$. Let $h = cn^{-1/7}$ $(c > 0)$ be the bandwidth used to compute $\hat{\lambda}_n^{\mathrm{SM}}$. Then*

$$n^{2/7}(\hat{\lambda}_n^{\mathrm{SM}}(t) - \lambda_0(t)) \rightsquigarrow \mathcal{N}(\mu_{\lambda,\mathrm{SM}}, \sigma_{\lambda,\mathrm{SM}}^2),$$

*where*

$$\mu_{\lambda,\mathrm{SM}} = \frac{1/2c^2 m_2(k)}{1-F_0(t)}\left(f_0''(t) + \frac{f_0(t)f_0'(t)}{1-F_0(t)}\right),$$

$$\sigma^2_{\lambda,\mathrm{SM}} = \frac{F_0(t)}{c^3 g(t)(1-F_0(t))}\int k'(u)^2\,du$$

*for* $t$ *such that* $(1-F_0(t))f_0''(t) + f_0(t)f_0'(t) \neq 0$.

*For fixed* $t>0$ *the aMSE-optimal bandwidth* $h$ *for* $\hat\lambda_n^{\mathrm{SM}}(t)$ *is given by* $h_{n,\lambda,\mathrm{SM}} = \sigma^2_{\lambda,\mathrm{SM}} n^{-1/7}$, *where*

$$\begin{aligned} c_{\lambda,\mathrm{SM}} &= \left\{3\frac{F_0(t)}{g(t)(1-F_0(t))}\int k'(u)^2\,du\right\}^{1/7} \\ &\quad\times\left\{\frac{m_2^2(k)}{(1-F_0(t))^2}\left(f_0''(t) + \frac{f_0(t)f_0'(t)}{1-F_0(t)}\right)^2\right\}^{-1/7}. \end{aligned} \tag{4.7}$$

PROOF. The proof uses the same decomposition as the proof of Corollary 3.7, but now $n^{2/7}R_n(t) \xrightarrow{\mathcal{P}} \frac{1}{2}c^2 m_2(k) f_0'(t)$. This gives that

$$T_n(t) \xrightarrow{\mathcal{P}} \frac{1}{2}c^2 m_2(k) f_0(t)\frac{f_0'(t)}{(1-F_0(t))^2} = \mu_{F,\mathrm{SM}}\frac{f_0'(t)}{(1-F_0(t))^2}$$

and $\mu_{\lambda,\mathrm{SM}} = \mu_{f,\mathrm{SM}}/(1-F_0(t)) + \mu_{F,\mathrm{SM}} f_0(t)/(1-F_0(t))^2$. $\square$

**5. Bandwidth selection in practice.** In the previous sections, we derived the optimal bandwidths to estimate $\theta_0(F)$ [the unknown distribution function $F_0$, its density $f_0$ or the hazard rate $\lambda_0 = f_0/(1-F_0)$ at a point $t$] using two different smoothing methods. These optimal bandwidths can be written as $h_{n,\hat\theta(F)} = c_{\hat\theta(F)} n^{-\alpha}$ for some $\alpha > 0$ (either $1/5$ or $1/7$), where $c_{\hat\theta(F)}$ is defined as the minimizer of aMSE$(c)$ over all positive $c$. For example $\theta_0(F) = F_0(t)$ and $\hat\theta(F) = \hat F_n^{\mathrm{SM}}(t)$. However, the asymptotic mean squared error depends on the unknown distribution $F_0$, so $c_{\hat\theta(F)}$ and $h_{n,\hat\theta(F)}$ are unknown.

Several data dependent methods are known to overcome this problem by estimating the aMSE, e.g., the bootstrap method of Efron (1979) or plug-in methods where the unknown quantities, like $f_0$ or $f_0''$, in the aMSE are replaced by estimates [see, e.g., Sheather (1983)]. We use the smoothed bootstrap method, which is commonly used to estimate the bandwidth in density-type problems, see, e.g., Hazelton (1996) and González-Manteiga, Cao and Marron (1996).

For $\hat\theta(F) = \hat F_n^{\mathrm{MS}}(t)$ the smoothed bootstrap works as follows. Let $n$ be the sample size and $h_0 = c_0 n^{-1/5}$ an initial choice of the bandwidth. Instead of sampling from the empirical distribution (as is done in the usual bootstrap) we sample $X_1^{*,1}, X_2^{*,1}, \ldots, X_m^{*,1}$ $(m \le n)$ from the distribution $\hat{F}_{n,h_0}^{\mathrm{SM}}$ (where we explicitly denote the bandwidth $h_0$ used to compute $\hat{F}_n^{\mathrm{SM}}$). Furthermore, we sample $T_1^{*,1}, \ldots, T_m^{*,1}$ from $\hat{G}_{n,h_0}$ and define $\Delta_i^{*,1} = 1_{\{X_i^{*,1} \le Z_i^{*,1}\}}$. Based on the sample $(T_1^{*,1}, \Delta_1^{*,1}), \ldots, (T_m^{*,1}, \Delta_m^{*,1})$, we determine the estimator $\hat{F}_{m,cm^{-1/5}}^{\mathrm{SM},1}$ with bandwidth $h = cm^{-1/5}$. We repeat this many times (say $B$ times), and estimate aMSE$(c)$ by

$$\widehat{\mathrm{MSE}}_B(c) = B^{-1} \sum_{i=1}^{B} (\hat{F}_{m,cm^{-1/5}}^{\mathrm{SM},i}(t) - \hat{F}_{n,h_0}(t))^2.$$

The optimal bandwidth $h_{n,F,\mathrm{SM}}$ we estimate by $\hat{h}_{n,F,\mathrm{SM}} = \hat{c}_{F,\mathrm{SM}} n^{-1/5}$ where $\hat{c}_{F,\mathrm{SM}}$ is defined as the minimizer of $\widehat{\mathrm{MSE}}_B(c)$ over all positive $c$. For the other estimators, the smoothed bootstrap works similarly.

Table 1 contains the values of $\hat{c}_{F,\mathrm{SM}}$ and $\hat{h}_{n,F,\mathrm{SM}}$ for the different choices of $c_0$ and two different points $t$ based on a simulation study. For the distribution of the $X_i$, we took a shifted Gamma(4) distribution, i.e., $f_0(x) = \frac{(x-2)^3}{3!} \exp(-(x-2)) 1_{[2,\infty)}(x)$, and for the distribution of the $T_i$ we took an exponential distribution with mean 3, i.e., $g(t) = \frac{1}{3}\exp(-t/3) 1_{[0,\infty)}$. Furthermore, we took $n = 10{,}000$, $m = 2000$, $B = 500$ and $k(t) = \frac{35}{32}(1-t^2)^3 1_{[-1,1]}(t)$, the triweight kernel. The table also contains the theoretical aMSE optimal values $c_{F,\mathrm{SM}}$, given in (4.5), the values of $\tilde{c}_{F,\mathrm{SM}}$ using Monte Carlo simulations of size $n = 10{,}000$ and $m = 2000$ and the corresponding values of $h_{n,F,\mathrm{SM}}$ and $\tilde{h}_{n,F,\mathrm{SM}}$. In the Monte Carlo simulation, we resampled $B$ times a sample of size $n$ (and $m$) from the true underlying distributions and estimated, in case of sample size $n$, the aMSE by

$$\widetilde{\mathrm{MSE}}_B(c) = B^{-1} \sum_{i=1}^{B} (\hat{F}_{n,cn^{-1/5},i}^{\mathrm{SM}}(t) - F_0(t))^2.$$

Then $\tilde{c}_{F,\mathrm{SM}}$ is defined as the minimizer of $\widetilde{\mathrm{MSE}}_B(c)$ over all positive $c$ and $\tilde{h}_{F,\mathrm{SM}} = \tilde{c}_{F,\mathrm{SM}} n^{-1/5}$. Figure 2 shows the aMSE$(c)$ for $t = 4$ and its estimates $\widehat{\mathrm{MSE}}_B(c)$ with $c_0 = 15$ and $\widetilde{\mathrm{MSE}}_B(c)$. Figure 2 also shows the estimator $\hat{F}_n^{\mathrm{SM}}$ with bandwidth $h = 1.7$ (which is somewhere in the middle of the results in Table 1 for $c_0 = 15$), the maximum likelihood estimator $\hat{F}_n$ and the true distribution $F_0$.

We also applied the smoothed bootstrap to choose the smoothing parameter for $\hat{F}_n^{\mathrm{SM}}(t)$ based on the hepatitis A prevalence data described by Keiding (1991). Table 2 contains the values of $\hat{c}_{F,\mathrm{SM}}$ and $\hat{h}_{n,F,\mathrm{SM}}$ for three different time points, $t = 20$, $t = 45$ and $t = 70$ and for different values of $c_0$.

TABLE 1
*Minimizing values for c and corresponding values of the bandwidth based on the smoothed bootstrap method for different values of $c_0$, based on Monte Carlo simulations and the theoretical values*

| | $t = 4.0$ | | $t = 6.5$ | |
|---|---|---|---|---|
| | $\hat{c}_{F,\mathrm{SM}}$ | $\hat{h}_{n,F,\mathrm{SM}}$ | $\hat{c}_{F,\mathrm{SM}}$ | $\hat{h}_{n,F,\mathrm{SM}}$ |
| $c_0 = 5$ | 6.050 | 0.959 | 9.150 | 1.450 |
| $c_0 = 10$ | 7.350 | 1.165 | 10.100 | 1.601 |
| $c_0 = 15$ | 7.700 | 1.220 | 12.050 | 1.910 |
| $c_0 = 20$ | 7.850 | 1.244 | 14.150 | 2.243 |
| $c_0 = 25$ | 9.850 | 1.561 | 15.500 | 2.457 |
| MC-sim ($n$) | 6.700 | 1.062 | 10.700 | 1.696 |
| MC-sim ($m$) | 6.750 | 1.070 | 11.600 | 1.838 |
| Theor. val. | 6.467 | 1.025 | 10.426 | 1.652 |

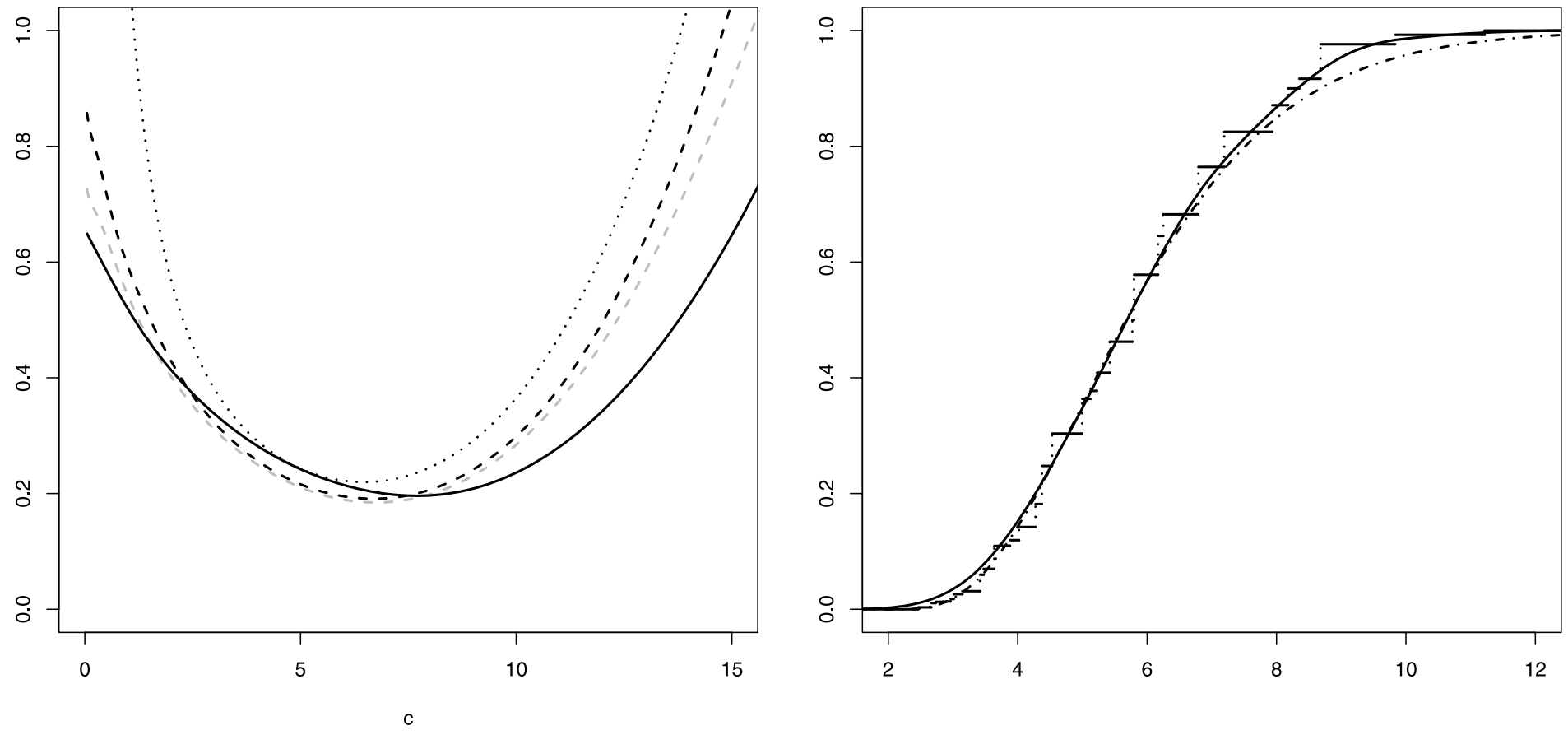

FIG. 2. *Left panel: the aMSE of $\hat{F}_n^{\mathrm{SM}}(4)$ (dotted line) and its estimates based on the smoothed bootstrap (solid line) with $c_0 = 15$ and the Monte Carlo simulations (dashed lines) with sample size $n$ (black line) and $m$ (grey line). Right panel: the true distribution (dash-dotted line) and its estimators $\hat{F}_n^{\mathrm{SM}}$ with $h = 1.7$ (solid line) and $\hat{F}_n$ (step function).*

The size $n$ of the hepatitis A prevalence data is 850. For the sample size $m$ of the smoothed bootstrap sample, we took 425 and we repeated the smoothed bootstrap $B = 500$ times. If we take the smoothing parameter $h$ equal to 25 (which is somewhere in the middle of the results in Table 2), the resulting estimator $\hat{F}_n^{\mathrm{SM}}$ is shown in Figure 3. The maximum likelihood estimator $\hat{F}_n$ is also shown in Figure 3.

TABLE 2

*Minimizing values for c and corresponding values of the bandwidth based on the smoothed bootstrap method for different values of $c_0$ and for three different values of t*

| | $t = 20$ | | $t = 45$ | | $t = 70$ | |
|---|---|---|---|---|---|---|
| | $\hat{c}_{F,\mathrm{SM}}$ | $\hat{h}_{n,F,\mathrm{SM}}$ | $\hat{c}_{F,\mathrm{SM}}$ | $\hat{h}_{n,F,\mathrm{SM}}$ | $\hat{c}_{F,\mathrm{SM}}$ | $\hat{h}_{n,F,\mathrm{SM}}$ |
| $c_0 = 50$ | 107.7 | 27.947 | 60.3 | 15.647 | 128.9 | 33.448 |
| $c_0 = 60$ | 105.6 | 27.402 | 67.6 | 17.541 | 128.7 | 33.396 |
| $c_0 = 70$ | 106.7 | 27.687 | 67.8 | 17.593 | 127.4 | 33.059 |
| $c_0 = 80$ | 101.8 | 26.416 | 71.6 | 18.579 | 130.4 | 33.837 |
| $c_0 = 90$ | 92.5 | 24.003 | 70.4 | 18.268 | 131.0 | 33.993 |
| $c_0 = 100$ | 91.9 | 23.847 | 76.5 | 19.851 | 127.5 | 33.085 |
| $c_0 = 110$ | 90.5 | 23.484 | 75.9 | 19.695 | 126.2 | 32.747 |
| $c_0 = 120$ | 89.8 | 23.302 | 80.8 | 20.967 | 124.3 | 32.254 |
| $c_0 = 130$ | 89.4 | 23.198 | 81.0 | 21.018 | 124.5 | 32.306 |
| $c_0 = 140$ | 84.2 | 21.849 | 81.9 | 21.252 | 120.2 | 31.190 |
| $c_0 = 150$ | 87.3 | 22.653 | 88.7 | 23.017 | 117.4 | 30.464 |

**6. Discussion.** We considered two different methods to obtain smooth estimates for the distribution function $F_0$ and its density $f_0$ in the current status model. Pointwise asymptotic results show that for estimating any of these functions both estimators have the same variance but a different asymptotic bias. The asymptotic bias of the MSLE equals the asymptotic bias of the SMLE plus an additional term depending on the unknown densities $f_0$ and $g$ (and their derivatives) and the point $t$ we estimate at. For some choices of $f_0$ and $g$ this additional term is positive, for other choices it is negative. Hence, we cannot say one method always results in a smaller bias than the other method, i.e., one estimator is uniformly superior. This was also seen by Marron and Padgett (1987) and Patil, Wells and Marron (1994) in the case of estimating densities based on right-censored data. Figure 4 shows the asymptotic mean squared error of the estimators $\hat{F}_n^{\mathrm{MS}}(t)$ and $\hat{F}_n^{\mathrm{SM}}(t)$ if $F_0$ is the shifted Gamma(4) distribution and $G$ is the exponential distribution with mean 3, i.e., $f_0(x) = \frac{(x-2)^3}{3!}\exp(-(x-2))1_{[2,\infty)}(x)$, $g(t) = \frac{1}{3}\exp(-t/3)1_{[0,\infty)}$ and $c = 7.5$. For some values of $t$ the aMSE of $\hat{F}_n^{\mathrm{MS}}(t)$ is smaller [meaning that the bias of $\hat{F}_n^{\mathrm{MS}}(t)$ is smaller], for other values of $t$ the aMSE of $\hat{F}_n^{\mathrm{SM}}(t)$ is smaller [meaning that the bias of $\hat{F}_n^{\mathrm{SM}}(t)$ is smaller].

We also considered smooth estimators for the hazard rate $\lambda_0$, defined as

$$\hat{\lambda}_n(t) = \frac{\hat{f}_n(t)}{1 - \hat{F}_n(t)},$$

where $\hat{f}_n$ and $\hat{F}_n$ are either $\hat{f}_n^{\mathrm{MS}}$ and $\hat{F}_n^{\mathrm{MS}}$ or $\hat{f}_n^{\mathrm{SM}}$ and $\hat{F}_n^{\mathrm{SM}}$. Because $\hat{\lambda}_n(t)$ is a quotient, we could estimate nominator and denominator separately by choosing one bandwidth $h = cn^{-1/7}$ to compute $\hat{f}_n(t)$ and a different bandwidth $h_1 = c_1 n^{-1/5}$ to compute $\hat{F}_n(t)$. However, by the relation

$$\lambda_0(t) = \frac{d}{dz} - \log(1 - F_0(z))|_{z=t} = \frac{f_0(t)}{1 - F_0(t)}$$

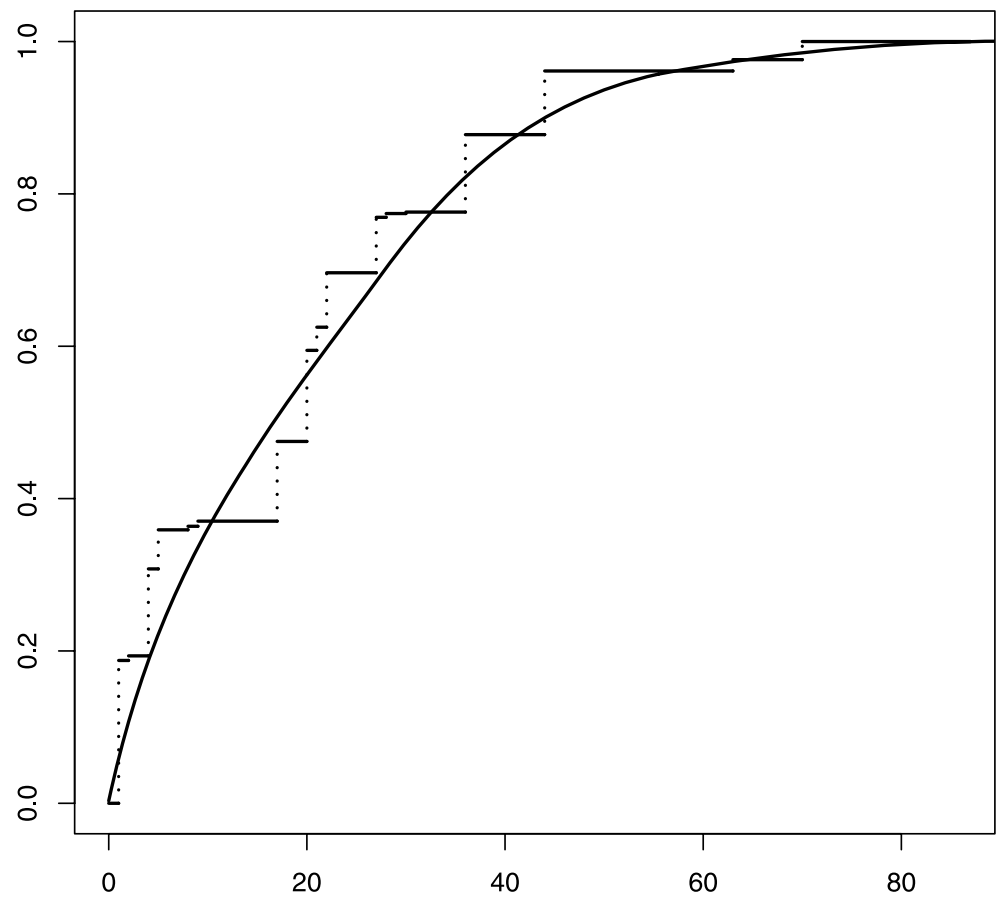

FIG. 3. *The estimators $\hat{F}_n^{\mathrm{SM}}$ (solid line) and $\hat{F}_n$ (dashed line) for the hepatitis A prevalence data.*

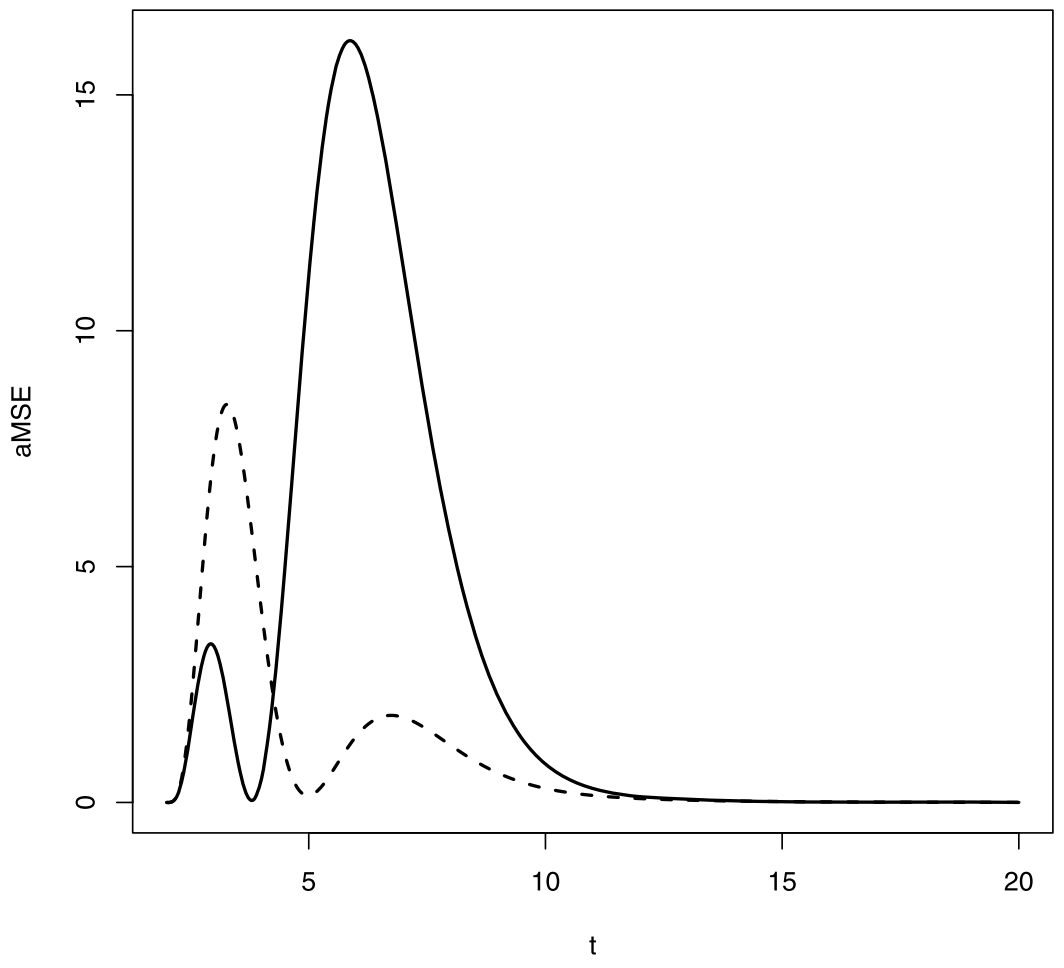

FIG. 4. *The aMSE of $\hat{F}_n^{\mathrm{MS}}(t)$ (solid line) and $\hat{F}_n^{\mathrm{SM}}(t)$ (dashed line) as function of $t$ in the situation described in Section 6.*

it is more natural to estimate $f_0(t)$ and $F_0(t)$ with the same bandwidth. As for the estimators for $f_0$ and $F_0$, we cannot say the estimator $\hat{\lambda}_n^{\mathrm{MS}}(t)$ with bandwidth of order $n^{-1/7}$ is uniformly superior to $\hat{\lambda}_n^{\mathrm{SM}}(t)$ with bandwidth of order $n^{-1/7}$.

## APPENDIX: TECHNICAL LEMMAS AND PROOFS

In this section, we prove most of the results stated in the previous sections. We start with some results on the consistency and pointwise asymptotics of the kernel estimators $\hat{g}_n$, $\hat{g}_n'$, $\hat{G}_n$, $\hat{g}_{n,1}$, $\hat{g}_{n,1}'$ and $\hat{G}_{n,1}$.

LEMMA A.1. *Let $\hat{g}_n$ be the boundary kernel estimator for $g$, with smoothing parameter $h = n^{-\alpha}$ ($\alpha < 1/3$). Then with probability converging to one $\hat{g}_n$ is uniformly bounded, i.e.,*

$$\exists C > 0 : P\Bigl(\sup_{x\in[0,1]} |\hat{g}_n(x)| \le C\Bigr) \longrightarrow 1. \tag{A.1}$$

PROOF. First note that without loss of generality we can assume $0 \le k(u) \le k(0)$. Recall that $\nu_{i,\beta}(k) = \int_{-1}^{\beta} u^i k(u)\,du$ for $\beta \in [0,1]$, for which we have the following bounds

$$\nu_{0,\beta} \ge \tfrac{1}{2}, \qquad |\nu_{1,\beta}| \le \tfrac{1}{2}\mathrm{E}_k|U|, \qquad \tfrac{1}{2}\operatorname{Var}_k U \le \nu_{2,\beta} \le \operatorname{Var}_k U,$$

where $U$ has density $k$. Combining this, we get that $\nu_{0,\beta}\nu_{2,\beta} - \nu_{1,\beta}^2 \ge \frac{1}{4}\operatorname{Var}_k |U| > 0$, so that we can uniformly bound the kernel $k^\beta$ by

$$\begin{aligned}
|k^\beta(u)| &= \left|\frac{\nu_{2,\beta} - \nu_{1,\beta}u}{\nu_{0,\beta}\nu_{2,\beta} - \nu_{1,\beta}^2} k(u) 1_{(-1,\beta]}(u)\right| \\
&\le \frac{|\nu_{2,\beta}| + |\nu_{1,\beta}|}{1/4\operatorname{Var}_k |U|} k(0) = ck(0).
\end{aligned}$$

For the boundary kernel estimate $\hat{g}_n$, we then have

$$\begin{aligned}
|\hat{g}_n(x)| &= \left|h^{-1}\int k^\beta((x-y)/h)\,d\mathbb{G}_n(y)\right| \\
&\le h^{-1}ck(0)|\mathbb{G}_n(x+h) - \mathbb{G}_n(x-h)| \\
&\le h^{-1}ck(0)|\mathbb{G}_n(x+h) - G(x+h) - \mathbb{G}_n(x-h) + G(x-h)| \\
&\quad + h^{-1}ck(0)(G(x+h) - G(x-h)) \\
&\le ck(0)n^{\alpha-1/2} 2\sup_{y\ge 0}\sqrt{n}|\mathbb{G}_n(y) - G(y)| + 2\|g\|_\infty ck(0) \\
&= \mathcal{O}_p(n^{\alpha-1/2}) + 2\|g\|_\infty ck(0).
\end{aligned}$$

Since this bound in uniform in $x$, (A.1) follows for $C = 3\|g\|_\infty ck(0)$. □

Lemma A.2. *Assume $g$ satisfies conditions* (G.1) *and let* $\hat{G}_n, \hat{G}_{n,1}, \hat{g}_n, \hat{g}_{n,1}$, $\hat{g}_n'$ *and* $\hat{g}_{n,1}'$ *be kernel estimators for* $G$, $G_1$, $g$, $g_1$, $g'$ *and* $g_1'$ *with kernel density* $k$ *satisfying condition* (K.1) *and bandwidth* $h = cn^{-\alpha}$ $(c > 0)$. *For* $\alpha \in (0, 1/3)$ *and* $m > 0$

$$\sup_{t\in[m,\infty)} |\hat{g}_n(t) - g(t)| \xrightarrow{\mathcal{P}} 0, \qquad \sup_{t\in[m,\infty)} |\hat{g}_n'(t) - g'(t)| \xrightarrow{\mathcal{P}} 0,$$
$$\sup_{t\in[0,2M_0]} |\hat{G}_n(t) - G(t)| \xrightarrow{\mathcal{P}} 0, \tag{A.2}$$

$$\sup_{t\in[m,\infty)} |\hat{g}_{n,1}(t) - g_1(t)| \xrightarrow{\mathcal{P}} 0, \qquad \sup_{t\in[m,\infty)} |\hat{g}_{n,1}'(t) - g_1'(t)| \xrightarrow{\mathcal{P}} 0,$$
$$\sup_{t\in[0,2M_0]} |\hat{G}_{n,1}(t) - G_1(t)| \xrightarrow{\mathcal{P}} 0. \tag{A.3}$$

Proof. Let $\hat{g}_n^u$ be the uncorrected kernel estimate for $g$ and note that by properties of the boundary kernel estimator we have for all $x \geq h$

$$\hat{g}_n^u(x) = \hat{g}_n(x).$$

Hence, the first two results in (A.2) follow immediately from Theorems A and C in Silverman (1978). To prove the third result in (A.2), fix $M > M_0$, $\epsilon > 0$ and choose $0 < \delta < \varepsilon/(2C)$ such that $G(\delta) < \varepsilon/4$, where $C$ is such that (A.1) holds. For all $x \geq 0$ and $n$ sufficiently large (such that $h = h_n < \delta$), we then have

$$|\hat{G}_n(x) - G(x)| \leq \delta \sup_{y\in[0,\delta]} |\hat{g}_n(y)| + G(\delta) + \sup_{y\geq\delta} |\hat{G}_n(y) - G(y)|.$$

The right-hand side does not depend on $x$ so that

$$\begin{aligned}
&P(\|\hat{G}_n - G\|_\infty > \epsilon) \\
&\quad\leq P\Big(\delta \sup_{y\in[0,\delta]} |\hat{g}_n(y)| + G(\delta) + \sup_{y\geq\delta} |\hat{G}_n^u(y) - G(y)| > \epsilon\Big) \\
&\quad\leq P\Big(\delta \sup_{y\in[0,1]} |\hat{g}_n(y)| + G(\delta) + \sup_{y\geq\delta} |\hat{G}_n^u(y) - G(y)| > \epsilon\Big) \\
&\quad= P\Big(\Big\{\delta \sup_{y\in[0,1]} |\hat{g}_n(y)| + G(\delta) + \sup_{y\geq\delta} |\hat{G}_n^u(y) - G(y)| > \epsilon\Big\} \\
&\qquad\qquad\qquad\qquad \cap \Big\{\sup_{y\in[0,1]} |\hat{g}_n(y)| \leq C\Big\}\Big) \\
&\qquad + P\Big(\Big\{\delta \sup_{y\in[0,1]} |\hat{g}_n(y)| + G(\delta) + \sup_{y\geq\delta} |\hat{G}_n^u(y) - G(y)| > \epsilon\Big\}
\end{aligned}$$

$$\cap \Big\{ \sup_{y\in[0,1]} |\hat g_n(y)| > C \Big\} \Big)$$

$$\le P\Big(\sup_{y\ge\delta} |\hat G_n^u(y) - G(y)| > \epsilon/4\Big).$$

The last probability converges to zero as a consequence of Theorem A in Silverman (1978), hence $\|\hat G_n - G\|_\infty \xrightarrow{\mathcal{P}} 0$.

For the first result in (A.3), define a binomially distributed random variable $N_1 = \sum_{i=1}^n \Delta_i$ with parameters $n$ and $p = P(\Delta_1 = 1) = \int F_0(u) g(u)\,du$, and the probability density $\tilde g(t) = g_1(t)/p$. Let $V_1, \dots, V_{N_1}$ be the $T_i$ such that $\Delta_i = 1$, and rewrite $\hat g_{n,1}(t)$ as $\frac{1}{nh}\sum_{i=1}^{N_1} k_h(t - V_i) = \frac{N_1}{n}\hat g_{N_1}(t)$. Then we have by the triangle inequality

$$\|\hat g_{n,1} - g_1\|_\infty = \left\| \hat g_{N_1}\frac{N_1}{n} - \tilde g p \right\|_\infty \le p\|\hat g_{N_1} - \tilde g\|_\infty + \|\hat g_{N_1}\|_\infty \left| \frac{N_1}{n} - p \right|.$$

The first term on the right-hand side converges to zero in probability by Silverman (1978), since $N_1 \xrightarrow{\mathcal{P}} \infty$ as $n \to \infty$. For the second term on the right-hand side, note that

$$\|\hat g_{N_1}\|_\infty = \|\tilde g + \hat g_{N_1} - \tilde g\|_\infty \le \|\tilde g\|_\infty + \|\hat g_{N_1} - \tilde g\|_\infty,$$

where the last term again converges to zero in probability by Silverman (1978). Combining this with the Law of Large Numbers applied to $|\frac{N_1}{n} - p|$ gives that $\|\hat g_{N_1}\|_\infty |\frac{N_1}{n} - p| \xrightarrow{\mathcal{P}} 0$ as $n \to \infty$, hence $\|\hat g_{n,1} - g_1\|_\infty \xrightarrow{\mathcal{P}} 0$. The proofs of the other results in (A.3) are similar. $\square$

LEMMA A.3. *Let $\hat g_n$ and $\hat g_{n,1}$ be kernel estimates for $g$ and $g_1$ with kernel density $k$ satisfying condition* (K.1) *and bandwidth $h = cn^{-\alpha}$ $(c > 0)$. Fix $t > 0$ such that $f_0''$ and $g''$ exist and are continuous at $t$. Then for $\alpha = 1/5$,*

$$n^{2/5}\left( \begin{pmatrix} \hat g_n(t) \\ \hat g_{n,1}(t) \end{pmatrix} - \begin{pmatrix} g(t) \\ g_1(t) \end{pmatrix} \right) \rightsquigarrow \mathcal{N}\left( \begin{pmatrix} \frac12 c^2 m_2(k) g''(t) \\ \frac12 c^2 m_2(k) g_1''(t) \end{pmatrix}, \Sigma_1 \right) \tag{A.4}$$

*with*

$$\Sigma_1 = c^{-1}\int k(u)^2\,du \begin{pmatrix} g(t) & g_1(t) \\ g_1(t) & g_1(t) \end{pmatrix}. \tag{A.5}$$

*For $0 < \alpha < 1/5$,*

$$n^{2\alpha}(\hat g_n(t) - g(t)) \xrightarrow{\mathcal{P}} \tfrac12 c^2 m_2(k) g''(t)$$

*and*

$$n^{2\alpha}(\hat g_{n,1}(t) - g_1(t)) \xrightarrow{\mathcal{P}} \tfrac12 c^2 m_2(k) g_1''(t).$$

*Let $\hat{g}'_{n,1}$ and $\hat{g}'_n$ be as defined in (3.7). Then for fixed $t>0$ such that $f_0^{(3)}$ and $g^{(3)}$ exist and are continuous at $t$ and $\alpha=1/7$,*

$$n^{2/7}\left(\begin{pmatrix}\hat{g}'_n(t)\\ \hat{g}'_{n,1}(t)\end{pmatrix}-\begin{pmatrix}g'(t)\\ g'_1(t)\end{pmatrix}\right)\rightsquigarrow\mathcal{N}\left(\begin{pmatrix}\frac{1}{2}c^2m_2(k)g^{(3)}(t)\\ \frac{1}{2}c^2m_2(k)g_1^{(3)}(t)\end{pmatrix},\Sigma_2\right) \tag{A.6}$$

*with*

$$\Sigma_2=c^{-3}\int k'(u)^2\,du\begin{pmatrix}g(t)&g_1(t)\\ g_1(t)&g_1(t)\end{pmatrix}. \tag{A.7}$$

Proof. We start with the proof of (A.4). Define

$$Y_i=\begin{pmatrix}Y_{i;1}\\ Y_{i;2}\end{pmatrix}=n^{-3/5}\begin{pmatrix}k_h(t-T_i)\\ k_h(t-T_i)\Delta_i\end{pmatrix}.$$

By the assumptions on $f_0$ and $g$ and condition (K.1), we have

$$\mathrm{E}Y_i=n^{-3/5}\begin{pmatrix}g(t)+\frac{1}{2}h^2m_2(k)g''(t)+\mathcal{O}_p(h^2)\\ g_1(t)+\frac{1}{2}h^2m_2(k)g''_1(t)+\mathcal{O}_p(h^2)\end{pmatrix},$$

$$\sum_{i=1}^n\operatorname{Var}Y_i=c^{-1}\int k(u)^2\,du\begin{pmatrix}g(t)&g_1(t)\\ g_1(t)&g_1(t)\end{pmatrix}+\mathcal{O}_p(n^{-1/5}).$$

By the Lindeberg–Feller central limit theorem, we get

$$n^{2/5}\left(\begin{pmatrix}\hat{g}_n(t)\\ \hat{g}_{n,1}(t)\end{pmatrix}-\begin{pmatrix}g(t)\\ g_1(t)\end{pmatrix}\right)-\begin{pmatrix}\frac{1}{2}c^2m_2(k)g''(t)\\ \frac{1}{2}c^2m_2(k)g''_1(t)\end{pmatrix}\rightsquigarrow\mathcal{N}(0,\Sigma_1),$$

where $\Sigma_1$ is defined in (A.5).

To prove that $n^{2\alpha}(\hat{g}_n(t)-g(t))\xrightarrow{\mathcal{P}}\frac{1}{2}c^2m_2(k)g''(t)$ for $0<\alpha<1/5$, define $W_i=n^{2\alpha-1}k_h(t-T_i)$. Since we have

$$\mathrm{E}W_i=n^{2\alpha-1}(g(t)+\tfrac{1}{2}h^2g''(t)+\mathcal{O}_p(h^2)),$$

$$n\operatorname{Var}W_i=n^{5\alpha-1}c^{-1}g(t)\int k(u)^2\,du+\mathcal{O}_p(n^{4\alpha-1})=\mathcal{O}_p(n^{5\alpha-1}),$$

we have that $\sum\operatorname{Var}W_i\longrightarrow 0$ for $0<\alpha<1/5$, hence

$$n\left(\frac{1}{n}\sum_{i=1}^n W_i-\mathrm{E}W_1\right)=n^{2\alpha}(\hat{g}_n(t)-g(t))-\frac{1}{2}c^2m_2(k)g''(t)+\mathcal{O}_p(1)\xrightarrow{\mathcal{P}}0.$$

Similarly we can prove that $n^{2\alpha}(\hat{g}_{n,1}(t)-g_1(t))\xrightarrow{\mathcal{P}}\frac{1}{2}c^2m_2(k)g''_1(t)$.

The proof of (A.6) is similar as the proof of (A.4). □

Using these results we now can prove the results in Section 3.

Proof of Lemma 3.2. The proof of the inequalities in Lemma 3.2 is based on the Monotone Convergence theorem (MCT). Denote the lower convex hull of the continuous cusum diagram defined in (3.2) by $t \mapsto (\hat{G}_n(t), C_n(t))$ for $t \in [0, \tau]$, where $\tau = \sup\{t \geq 0 : \hat{g}_{n,0}(t) + \hat{g}_{n,1}(t) > 0\}$. By definition of this convex hull, we have for all $t > 0$

$$\begin{aligned}\hat{G}_{n,1}(t) &= \int 1_{[0,t]}(u)\, d\hat{G}_{n,1}(u) \geq \int 1_{[0,t]}(u)\, dC_n(u) \\ &= \int \hat{F}_n^{\mathrm{MS}}(u) 1_{[0,t]}(u)\, d\hat{G}_n(u).\end{aligned} \tag{A.8}$$

The function $1_{[0,t]}(u)$ is decreasing on $[0, \infty)$. Consider an arbitrary distribution function $F$ on $[0, \infty)$ and write $p(t) = -\log F(t)$. Then, on $[0, \tau]$, the function $p$ can be approximated by decreasing step functions

$$p_m(t) = \sum_{i=1}^m a_i 1_{[0,x_i]}(t) \qquad \text{with } a_i \geq 0\ \forall i \text{ and } 0 < x_1 < \cdots < x_m < \tau.$$

The functions $p_m$ can be taken such that $p_m \uparrow p$, on $[0, \tau]$. For each $m$, we have

$$\begin{aligned}\int p_m(t)\, d\hat{G}_{n,1}(t) &= \sum_{i=1}^m \int a_i 1_{[0,x_i]}(t)\, d\hat{G}_{n,1}(t) \\ &\geq \sum_{i=1}^m \int a_i 1_{[0,x_i]}(t)\, dC_n(t) \\ &= \int p_m(t) \hat{F}_n^{\mathrm{MS}}(t)\, d\hat{G}_n(t).\end{aligned} \tag{A.9}$$

The MCT now gives that for each $n$

$$\begin{aligned}\lim_{m\to\infty} \int p_m(t)\, d\hat{G}_{n,1}(t) &= \int p(t)\, d\hat{G}_{n,1}(t) = -\int \log F(t)\, d\hat{G}_{n,1}(t), \\ \lim_{m\to\infty} \int p_m(t)\, dC_n(t) &= \int p(t)\, dC_n(t) = -\int \hat{F}_n^{\mathrm{MS}}(t) \log F(t)\, d\hat{G}_n(t).\end{aligned}$$

Combined with (A.9), this implies the first inequality in Lemma 3.2.

To prove the second inequality in Lemma 3.2, it suffices to prove

$$\int \log(1 - F(t))\, d\hat{G}_{n,1}(t) \geq \int \hat{F}_n^{\mathrm{MS}}(t) \log(1 - F(t))\, d\hat{G}_n(t), \tag{A.10}$$

since

$$\int \log(1 - F(t))\, d\hat{G}_{n,0}(t) = \int \log(1 - F(t))\, d(\hat{G}_n - \hat{G}_{n,1})(t).$$

The proof of (A.10) follows by a similar argument. Then we use approximations $q_m(t)$ of the decreasing function $q(t) = \log(1 - F(t))$ such that $q_m \uparrow q$ to prove (A.10).

For the equality statements for $F = \hat{F}_n^{\mathrm{MS}}$ in Lemma 3.2, we can also use the monotone approximation by step functions, restricting the jumps to the points of increase of $\hat{F}_n^{\mathrm{MS}}$ [i.e., points $x$ for which $\hat{F}_n^{\mathrm{MS}}(x+\epsilon) - \hat{F}_n^{\mathrm{MS}}(x-\epsilon) > 0$ for all $\epsilon > 0$] implying equality in (A.9). $\square$

Proof of Theorem 3.3. Take $0 < m < M < M_0$. By assumption (G.1) and Lemma A.2, with probability arbitrarily close to one, we have for $n$ sufficiently large that $\hat{g}_n(t) > 0$ for all $t \in [m, M]$. We then have that $\hat{F}_n^{\mathrm{naive}}(t) = \hat{g}_{n,1}(t)/\hat{g}_n(t)$ is well defined on $[m, M]$ and to prove that $\hat{F}_n^{\mathrm{naive}}(t)$ is monotonically increasing on $[m, M]$ with probability tending to one, it suffices to show that $\exists \delta > 0$ such that $\forall \eta > 0$

$$P\left(\forall t \in [m, M] : \frac{d}{dt}\hat{F}_n^{\mathrm{naive}}(t) \geq \delta\right) \geq 1 - \eta \tag{A.11}$$

for $n$ sufficiently large. We have that

$$\frac{d}{dt}\hat{F}_n^{\mathrm{naive}}(t) = \frac{\hat{g}_n(t)\hat{g}_{n,1}'(t) - \hat{g}_{n,1}(t)\hat{g}_n'(t)}{[\hat{g}_n(t)]^2},$$

which is also well defined.

To prove (A.11) it suffices to prove $\exists \delta > 0$ such that $\forall \eta > 0$

$$P(\forall t \in [m, M] : \hat{g}_n(t)\hat{g}_{n,1}'(t) - \hat{g}_{n,1}(t)\hat{g}_n'(t) \geq \delta) \geq 1 - \eta \tag{A.12}$$

for $n$ sufficiently large. For this, we write

$$\begin{aligned}
&\hat{g}_n(t)\hat{g}_{n,1}'(t) - \hat{g}_{n,1}(t)\hat{g}_n'(t) \\
&\quad = \hat{g}_n(t)(\hat{g}_{n,1}'(t) - g_1'(t)) + \hat{g}_{n,1}(t)(g'(t) - \hat{g}_n'(t)) \\
&\qquad + g_1'(t)(\hat{g}_n(t) - g(t)) + g'(t)(g_1(t) - \hat{g}_{n,1}(t)) + g(t)g_1'(t) - g'(t)g_1(t) \\
&\quad \geq - \sup_{t\in[m,M]} |\hat{g}_{n,1}'(t) - g_1'(t)| \sup_{t\in[m,M]} \hat{g}_n(t) \\
&\qquad - \sup_{t\in[m,M]} |\hat{g}_n'(t) - g'(t)| \sup_{t\in[m,M]} \hat{g}_{n,1}(t) \\
&\qquad - \sup_{t\in[m,M]} |\hat{g}_n(t) - g(t)| \sup_{t\in[m,M]} g_1'(t) \\
&\qquad - \sup_{t\in[m,M]} |\hat{g}_{n,1}(t) - g_1(t)| \sup_{t\in[m,M]} g'(t) \\
&\qquad + g^2(t)f_0(t).
\end{aligned}$$

By Lemma A.2 and assumptions (F.1) and (G.1), we have that (A.12) follows for $\delta < \inf_{t\in[m,M]} g^2(t)f_0(t)$. $\square$

Proof of Corollary 3.4. Fix $\delta > 0$ arbitrarily. We will prove that for $n$ sufficiently large

$$P(\hat{F}_n^{\text{naive}}(t) = \hat{F}_n^{\text{MS}}(t) \text{ for all } t \in [m, M]) \geq 1 - \delta.$$

Define for $\eta_1 \in (0, m)$, $\eta_2 \in (0, M_0 - M)$ and $n \geq 1$ the event $A_n$ by

$$A_n = \{\hat{F}_n^{\text{naive}}(t) \text{ is monotonically increasing and } \hat{g}_n(t) > 0 \\ \text{for } t \in [m - \eta_1, M + \eta_2]\}.$$

By Lemma A.2 and Theorem 3.3, we have for all $n$ sufficiently large $P(A_n) \geq 1 - \delta/10$.

Define the "linearly extended $\hat{G}_{n,1}$" by

$$C_n^*(t) = \begin{cases} \hat{G}_{n,1}(m) + (\hat{G}_n(t) - \hat{G}_n(m))\hat{F}_n^{\text{naive}}(m), & \text{for } t \in [0, m), \\ \hat{G}_{n,1}(t), & \text{for } t \in [m, M], \\ \hat{G}_{n,1}(M) + (\hat{G}_n(t) - \hat{G}_n(M))\hat{F}_n^{\text{naive}}(M), & \text{for } t \in (M, M_0]. \end{cases}$$

It now suffices to prove that for all $n$ sufficiently large

$$\text{(i)} \quad P(\{(\hat{G}_n(t), C_n^*(t)) : t \geq 0\} \text{ convex}) \geq 1 - \delta/2,$$

$$\text{(ii)} \quad P(\forall t \in [0, M_0] : C_n^*(t) \leq \hat{G}_{n,1}(t)) \geq 1 - \delta/2.$$

Indeed, then with probability $\geq 1 - \delta$ the curve $\{(\hat{G}_n(t), C_n^*(t)) : t \geq 0\}$ is *a* lower convex hull of the CCSD $\{(\hat{G}_n(t), \hat{G}_{n,1}(t)) : t \geq 0\}$ with $C_n^*(t) = \hat{G}_{n,1}(t)$ for all $t \in [m, M]$. From this, it follows that $C_n^*(t) = C_n(t)$ for all $t \in [m, M]$, hence also $C_n(t) = \hat{G}_{n,1}(t)$ for all $t \in [m, M]$. This implies that for $n$ sufficiently large

$$P\left(\forall t \in [m, M] : \hat{F}_n^{\text{naive}}(t) = \frac{d\hat{G}_{n,1}(t)}{d\hat{G}_n(t)} = \frac{dC_n(t)}{d\hat{G}_n(t)} = \hat{F}_n^{\text{MS}}(t)\right) \geq 1 - \delta.$$

We now prove (i). For the intervals $[0, m)$ and $(M, M_0]$ the curve $\{(\hat{G}_n(t), C_n^*(t)) : t \geq 0\}$ is the tangent line of the CCSD at the points $(\hat{G}_n(m), \hat{G}_{n,1}(m))$ and $(\hat{G}_n(M), \hat{G}_{n,1}(M))$, respectively, so on the event $A_n$ the curve is convex. This gives for $n$ sufficiently large

$$P(\{(\hat{G}_n(t), C_n^*(t)) : t \geq 0\} \text{ convex}) \geq P(A_n) \geq 1 - \delta/10 \geq 1 - \delta/2.$$

To prove (ii), we split up the interval $[0, M_0]$ in five different intervals $\mathcal{I}_1 = [0, m - \eta_1)$, $\mathcal{I}_2 = [m - \eta_1, m)$, $\mathcal{I}_3 = [m, M]$, $\mathcal{I}_4 = (M, M + \eta_2]$ and $\mathcal{I}_5 = (M + \eta_2, M_0]$ and prove that for $1 \leq i \leq 5$

$$P(C_i) = P(\forall t \in \mathcal{I}_i : C_n^*(t) \leq \hat{G}_{n,1}(t)) \geq 1 - \delta/10. \tag{A.13}$$

For $t \in \mathcal{I}_3$, $C_n^*(t) = \hat{G}_{n,1}(t)$, hence (A.13) holds trivially. For the interval $\mathcal{I}_2$, we use that

$$\hat{G}_{n,1}(u) - \hat{G}_{n,1}(v) = (\hat{G}_n(u) - \hat{G}_n(v))\hat{F}_n^{\text{naive}}(\xi) \tag{A.14}$$

for some $\xi \in [u, v]$ (depending on $u$ and $v$). This gives

$$\begin{aligned} &P(\forall t \in \mathcal{I}_2 : \hat{G}_{n,1}(t) - C_n^*(t) \geq 0) \\ &\quad = P(\forall t \in \mathcal{I}_2 : (\hat{G}_n(t) - \hat{G}_n(m))(\hat{F}_n^{\text{naive}}(\xi) - \hat{F}_n^{\text{naive}}(m)) \geq 0) \\ &\quad = P(\forall t \in \mathcal{I}_2 : \hat{F}_n^{\text{naive}}(\xi) - \hat{F}_n^{\text{naive}}(m) \leq 0) \geq P(A_n) \geq 1 - \delta/10. \end{aligned}$$

For $\mathcal{I}_4$, we can reason similarly.

Now consider (A.13) for $i = 1$. For every $t \in \mathcal{I}_1$, we have

$$\begin{aligned} &G_1(t) - G_1(m) - F_0(m)(G(t) - G(m)) \\ &\qquad = \int_t^m (F_0(m) - F_0(u))\,dG(u) \\ &\qquad \geq \int_{m-\eta_1}^m (F_0(m) - F_0(u))\,dG(u). \end{aligned}$$

This means we have

$$\begin{aligned} &\hat{G}_{n,1}(t) - C_n^*(t) \\ &\quad \geq \hat{G}_{n,1}(t) - G_1(t) + G_1(m) - \hat{G}_{n,1}(m) + F_0(m)(G(t) - \hat{G}_n(t)) \\ &\qquad + F_0(m)(\hat{G}_n(m) - G(m)) + (\hat{F}_n^{\text{naive}}(m) - F_0(m))(\hat{G}_n(m) - \hat{G}_n(t)) \\ &\qquad + \int_{m-\eta_1}^m (F_0(m) - F_0(u))\,dG(u) \\ &\quad \geq -2\|\hat{G}_{n,1} - G_1\|_\infty - 2\|\hat{G}_n - G\|_\infty - 2|\hat{F}_n^{\text{naive}}(m) - F_0(m)| \\ &\qquad + \int_{m-\eta_1}^m (F_0(m) - F_0(u))\,dG(u). \end{aligned}$$

By assumption (F.1), we have $\int_{m-\eta_1}^m (F_0(m) - F_0(u))\,dG(u) > 0$ so (A.13) follows for $i = 1$ by Lemma A.2 and the pointwise consistency of $\hat{F}_n^{\text{naive}}$.

For $i = 5$, the proof of (A.13) is similar as for $i = 1$. □

To prove the results in Section 4 and the results below, we use piecewise constant versions of the functions $\psi_{h,t}$ and $\varphi_{h,t}$ defined in (4.1). These functions are constant on the same intervals where the MLE $\hat{F}_n$ is constant. Denote these intervals by $J_i = [\tau_i, \tau_{i+1})$ for $0 \leq i \leq m - 1$ ($m \leq n$ and $\tau_0 = 0$) and the piecewise constant versions of $\psi_{h,t}$ and $\varphi_{h,t}$ by $\bar{\psi}_{h,t}$ and $\bar\varphi_{h,t}$. For $u \in J_i$ these functions can be written as $\bar\psi_{h,t}(u) = \psi(\hat A_n(u))$ and $\bar\varphi_{h,t}(u) = \varphi(\hat A_n(u))$ for $\hat A_n(u)$ defined as

$$\hat A_n(u) = \begin{cases} \tau_i, & \text{if } \forall t \in J_i : F_0(t) > \hat F_n(\tau_i), \\ s, & \text{if } \exists s \in J_i : \hat F_n(s) = F_0(s), \\ \tau_{i+1}, & \text{if } \forall t \in J_i : F_0(t) < \hat F_n(\tau_i), \end{cases} \tag{A.15}$$

for $u \in J_i$, see also Figure 5.

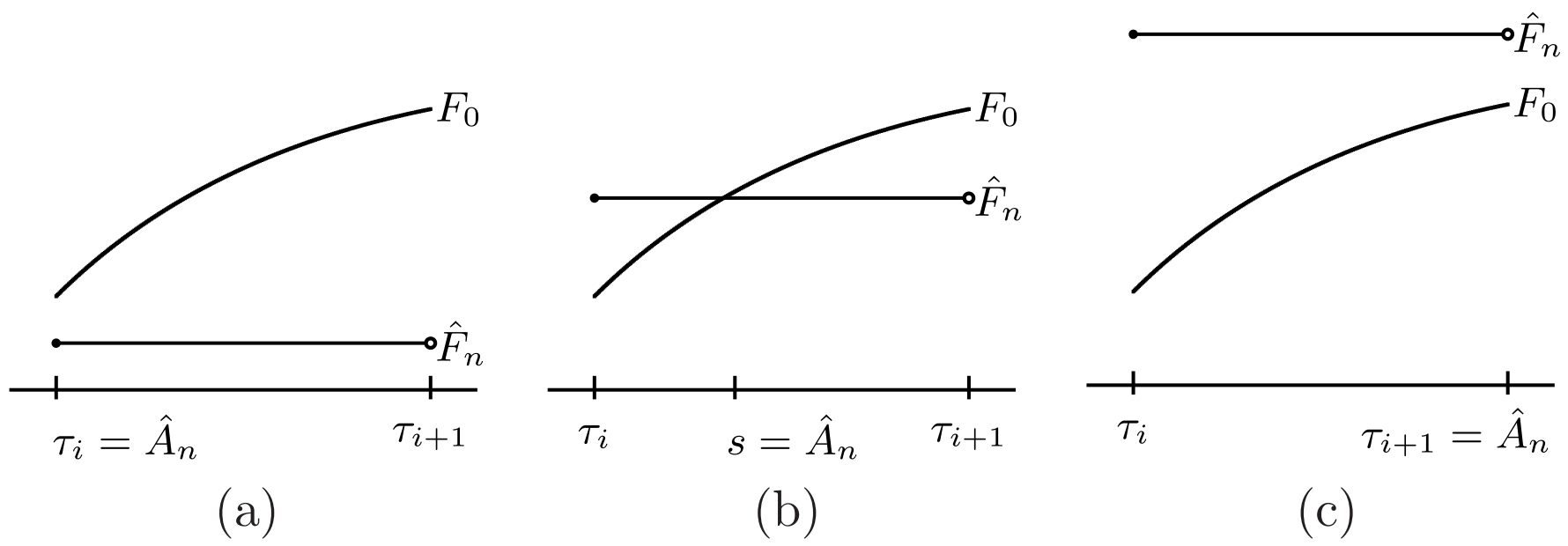

FIG. 5. *The 3 different possibilities for the function* $\hat A_n$. (a) $F_0(t) > \hat F_n(\tau_i)$; (b) $F_0(s) = \hat F_n(\tau_i)$; (c) $F_0(t) < \hat F_n(\tau_i)$.

We first derive upper bounds for the distance between the function $\psi_{h,t}$ and its piecewise constant version $\bar\psi_{h,t}$ and between $\varphi_{h,t}$ and $\bar\varphi_{h,t}$.

LEMMA A.4. *Let $t > 0$ be such that $f_0$ is positive and continuous in a neighborhood of $t$. Then there exists constants $c_1, c_2 > 0$ such that for $n$ sufficiently large*

$$|\bar\psi_{h,t}(u) - \psi_{h,t}(u)| \le \frac{c_1}{h^2}|\hat F_n(u) - F_0(u)| 1_{\{|t-u| \le h\}}, \tag{A.16}$$

$$|\bar\varphi_{h,t}(u) - \varphi_{h,t}(u)| \le \frac{c_2}{h^3}|\hat F_n(u) - F_0(u)| 1_{\{|t-u| \le h\}}. \tag{A.17}$$

PROOF. For $n$ sufficiently large, we have for all $s \in \mathcal{I}_t = [t-h, t+h]$ that $f_0(s) \ge \frac{1}{2} f_0(t)$. Fix $u \in \mathcal{I}_t$, then the interval $J_i$ it belongs to is of one of the following three types:

(i) $F_0(x) > \hat F_n(\tau_i)$ for all $x \in J_i$.
(ii) $F_0(x) = \hat F_n(x)$ for some $x \in J_i$.
(iii) $F_0(x) < \hat F_n(\tau_i)$ for all $x \in J_i$.

First, we consider the situation where $\hat F_n(u) = F_0(u)$. Then by definition of $\bar\psi_{h,t}$,

$$\bar\psi_{h,t}(u) = \psi_{h,t}(u),$$

so that both the left- and the right-hand side of (A.16) are equal to zero, and the upper bound holds. Note that for each $\hat{F}_n(u) = F_0(u)$ implies $\hat{A}_n(u) = u$, because $F_0$ is strictly increasing near $t$.

Now, we consider the situation where $\hat{F}_n(u) \neq F_0(u)$. For $v, \xi \in J_i$, we get by using a Taylor expansion

$$\begin{aligned} |\hat{F}_n(u) - F_0(u)| &= |\hat{F}_n(v) - F_0(u)| \\ &= |\hat{F}_n(v) - F_0(v) - (u - v) f_0(\xi)|. \end{aligned}$$

Now, we have three posibilities. If $\hat{A}_n(u) = \tau_i$, then we have that $F_0(\tau_i) - \hat{F}_n(\tau_i) \geq 0$ giving that

$$\begin{aligned} |\hat{F}_n(u) - F_0(u)| &= |\hat{F}_n(\tau_i) - F_0(\tau_i) - (u - \tau_i) f_0(\xi)| \\ &= |(u - \tau_i) f_0(\xi) + F_0(\tau_i) - \hat{F}_n(\tau_i)| \\ &\geq |u - \tau_i| f_0(\xi). \end{aligned}$$

If $\hat{A}_n(u) = v$ for some $v \neq u \in J_i$, then we have that $\hat{F}_n(v) = F_0(v)$, so that

$$\begin{aligned} |\hat{F}_n(u) - F_0(u)| = |\hat{F}_n(v) - F_0(u)| &= |\hat{F}_n(v) - F_0(v) - (u - v) f_0(\xi)| \\ &= |u - v| f_0(\xi). \end{aligned}$$

If $\hat{A}_n(u) = \tau_{i+1}$, then we have $\hat{F}_n(\tau_{i+1}-) - F_0(\tau_{i+1}) \geq 0$ giving that

$$\begin{aligned} |\hat{F}_n(u) - F_0(u)| &= |\hat{F}_n(\tau_{i+1}-) - F_0(\tau_{i+1}) - (u - \tau_{i+1}) f_0(\xi)| \\ &= |(\tau_{i+1} - u) f_0(\xi) + \hat{F}_n(\tau_{i+1}-) - F_0(\tau_{i+1})| \\ &\geq |\tau_{i+1} - u| f_0(\xi). \end{aligned}$$

For $v \in [\tau_i, \tau_{i+1}]$, this gives

$$|\hat{F}_n(u) - F_0(u)| \geq |u - v| f_0(\xi) \geq \tfrac{1}{2} f_0(t) |u - v| \geq 0.$$

Since it also holds that

$$\begin{aligned} |\bar{\psi}_{h,t}(u) - \psi_{h,t}(u)| &= |\psi_{h,t}(v) - \psi_{h,t}(u)| \leq c h^{-2} |v - u|, \\ |\bar{\varphi}_{h,t}(u) - \varphi_{h,t}(u)| &= |\varphi_{h,t}(v) - \varphi_{h,t}(u)| \leq \tilde{c} h^{-3} |v - u| \end{aligned}$$

the upper bound in (A.16) holds if $c_1 = 2c/f_0(t)$ and the upper bound in (A.17) holds if $c_2 = 2\tilde{c}/f_0(t)$. $\square$

To derive the asymptotic distribution of $\hat{F}_n^{\mathrm{SM}}(t)$ we need a result on the characterization of $\hat{F}_n$ and some results from empirical process theory, stated in Lemmas A.5 and A.7 below.

LEMMA A.5. *For every right continuous piecewise constant function $\varphi$ with only jumps at the points $\tau_1, \ldots, \tau_m$,*

$$\int \varphi(u)(\delta - \hat{F}_n(u))\, d\mathbb{P}_n(u, \delta) = 0.$$

PROOF. By the convex minorant interpretation of $\hat{F}_n$, we have that

$$\int_{[\tau_i, \tau_{i+1}) \times \{0,1\}} \delta\, d\mathbb{P}_n(u, \delta) = \int_{[\tau_i, \tau_{i+1}) \times \{0,1\}} \hat{F}_n(u)\, d\mathbb{P}_n(u, \delta)$$

for all $0 \le i \le m - 1$ (with $\tau_0 = 0$). This implies that

$$\begin{aligned}
&\int_{[\tau_i, \tau_{i+1}) \times \{0,1\}} \varphi(u)(\delta - \hat{F}_n(u))\, d\mathbb{P}_n(u, \delta) \\
&\qquad = \varphi(\tau_i) \int_{[\tau_i, \tau_{i+1}) \times \{0,1\}} (\delta - \hat{F}_n(u))\, d\mathbb{P}_n(u, \delta) = 0.
\end{aligned}$$

Hence,

$$\begin{aligned}
&\int \varphi(u)(\delta - \hat{F}_n(u))\, d\mathbb{P}_n(u, \delta) \\
&\qquad = \sum_{i=1}^{m-1} \int_{[\tau_i, \tau_{i+1}) \times \{0,1\}} \varphi(u)(\delta - \hat{F}_n(u))\, d\mathbb{P}_n(u, \delta) = 0.
\end{aligned}$$

□

Before we state the results on empirical process theory, we give some definitions and Theorem 2.14.1 in van der Vaart and Wellner (1996) needed for the proof of Lemma A.7.

Let $\mathcal{F}$ be the class of functions on $\mathbb{R}_+$ and $L_2(Q)$ the $L_2$-norm defined by a probability measure $Q$ on $\mathbb{R}_+$, i.e., for $g \in \mathcal{F}$

$$L_2(Q)[g] = \|g\|_{Q,2} = \left( \int_{\mathbb{R}_+} |g|\, dQ \right)^{1/2}.$$

For any probability measure $Q$, let $N(\varepsilon, \mathcal{F}, L_2(Q))$ be the minimal number of balls $\{g \in \mathcal{F} : \|g - f\|_{Q,2} < \varepsilon\}$ of radius $\varepsilon$ needed to cover the class $\mathcal{F}$. The entropy $H(\varepsilon, \mathcal{F}, L_2(Q))$ of $\mathcal{F}$ is then defined as

$$H(\varepsilon, \mathcal{F}, L_2(Q)) = \log N(\varepsilon, \mathcal{F}, L_2(Q))$$

and $J(\delta, \mathcal{F})$ is defined as

$$J(\delta, \mathcal{F}) = \sup_Q \int_0^\delta \sqrt{1 + H(\varepsilon, \mathcal{F}, L_2(Q))}\, d\varepsilon.$$

An envelope function of a function class $\mathcal{F}$ on $\mathbb{R}_+$ is any function $F$ such that $|f(x)| \le F(x)$ for all $x \in \mathbb{R}_+$ and $f \in \mathcal{F}$.

THEOREM A.6 [Theorem 2.14.1 in van der Vaart and Wellner (1996)]. *Let $P_0$ be the distribution of the observable vector $Z$ and $\mathcal{F}$ be a $P_0$-measurable class of measurable functions with measurable envelope function $F$. Then*

$$\mathrm{E}\sup_{f\in\mathcal{F}}\left|\int f\,d\sqrt{n}(\mathbb{P}_n-P_0)\right|\lesssim J(1,\mathcal{F})\|F\|_{P_0,2},$$

*where $\lesssim$ means $\le$ up to a multiplicative constant.*

LEMMA A.7. *Assume $F_0$ and $G$ satisfy conditions* (F.1) *and* (G.1) *and let $h:[0,\infty)\times\{0,1\}\to[-1,1]$ be defined as $h(u,\delta)=F_0(u)-\delta$. Then for $\alpha\le 1/5$ and $n\to\infty$*

$$R_n=n^{2\alpha}\int\bar{\psi}_{h,t}(u)(\hat{F}_n(u)-F_0(u))\,d(\mathbb{G}_n-G)(u)=\mathcal{O}_p(1), \tag{A.18}$$

$$S_n=n^{2\alpha}\int\{\bar{\psi}_{h,t}(u)-\psi_{h,t}(u)\}h(u,\delta)\,d(\mathbb{P}_n-P_0)(u,\delta)=\mathcal{O}_p(1). \tag{A.19}$$

PROOF. Define $\mathcal{I}_t=[t-\nu,t+\nu]$ for some $\nu>0$ and note that by (2.5) and (2.6) for any $\eta>0$ we can find $M_1,M_2>0$ such that for all $n$ sufficiently large

$$\begin{aligned}P(\mathcal{E}_{1,n,M_1})&:=P\Big(\sup_{u\in\mathcal{I}_t}|\hat{F}_n(u)-F_0(u)|\le M_1n^{-1/3}\log n\Big)\\&\ge 1-\eta/2,\end{aligned} \tag{A.20}$$

$$\begin{aligned}P(\mathcal{E}_{2,n,M_2})&:=P\Big(\sup_{u\in\mathcal{I}_t}|\hat{A}_n(u)-u|\le M_2n^{-1/3}\log n\Big)\\&\ge 1-\eta/2.\end{aligned} \tag{A.21}$$

Also note that $\|h\|_\infty\le 1$. Moreover, denote by $\mathcal{A}$ the class of monotone functions on $\mathcal{I}_t$, with values in $[0,2t]$. Then we know, see, e.g., (2.5) in van de Geer (2000), that for all $\delta>0$

$$H(\delta,\mathcal{A},L_2(Q))\lesssim\delta^{-1}$$

for any probability measure $Q$. For the same reason, the class $\mathcal{B}_M$ of functions of bounded variation on $[0,2t]$, absolutely bounded by $M$, has entropy function of the same order:

$$H(\delta,\mathcal{B}_M,L_2(Q))\lesssim\delta^{-1}\qquad\text{for all }\delta>0.$$

Let us now start the main argument. Choose $\eta>0$ and $M_1,M_2>0$ related to (A.20) and (A.21), correspondingly. Let $\nu_{1,n},\nu_{2,n}$ be vanishing sequences of positive numbers and write

$$\begin{aligned} P([|R_n| > \nu_{1,n}]) &= P([|R_n| > \nu_{1,n}] \cap \mathcal{E}_{1,n,M_1}) + P([|R_n| > \nu_{1,n}] \cap \mathcal{E}^c_{1,n,M_1}) \\ &\le P([|R_n| > \nu_{1,n}] \cap \mathcal{E}_{1,n,M_1}) + \eta/2 \le \nu_{1,n}^{-1} E|R_n| 1_{\mathcal{E}_{1,n,M_1}} + \eta/2, \end{aligned}$$

$$P([|S_n| > \nu_{2,n}]) \le P([|S_n| > \nu_{2,n}] \cap \mathcal{E}_{2,n,M_2}) + \eta/2 \le \nu_{2,n}^{-1} \mathrm{E}|S_n| 1_{\mathcal{E}_{2,n,M_2}} + \eta/2.$$

Here, we use the Markov inequality, (A.20) and (A.21). We now concentrate on the terms $\nu_{1,n}^{-1} E|R_n| 1_{\mathcal{E}_{1,n,M_1}}$ and $\nu_{2,n}^{-1}\mathrm{E}|S_n|1_{\mathcal{E}_{2,n,M_2}}$. We show that if we take, e.g., $\nu_{i,n} = \varepsilon n^{-\beta_i}(\log n)^2$ for $\beta_1 = 5/6 - 7\alpha/2$ and $\beta_2 = 5/6 - 4\alpha$ and any $\varepsilon > 0$ these terms will be smaller than $\eta/2$ for all $n$ sufficiently large, showing that $R_n = \mathcal{O}_p(n^{-\beta_1}(\log n)^2) = \mathcal{O}_p(1)$ and $S_n = \mathcal{O}_p(n^{-\beta_2}(\log n)^2) = \mathcal{O}_p(1)$ for $\alpha \le 1/5$.

We start with some definitions. Define for

$$C_n(u) = \frac{k(n^\alpha(t-u)/c)}{cg(u)} 1_{\mathcal{I}_t}(u),$$

the functions $\xi_{A,B,n}$ and $\zeta_{B,n}$ by

$$\begin{aligned} \xi_{A,B,n}(u) &= C_n(A(u))B(u), \\ \zeta_{B,n}(u,\delta) &= n^{1/3-\alpha}(\log n)^{-1} h(u,\delta)(C_n(n^{-1/3}B(u)\log n + u) - C_n(u)) \end{aligned}$$

and let

$$\mathcal{G}_{1,n} = \{\xi_{A,B,n} : A \in \mathcal{A}, B \in \mathcal{B}_{M_1}\}, \qquad \mathcal{G}_{2,n} = \{\zeta_{B,n} : B \in \mathcal{B}_{M_2}\}.$$

Note that by condition (K.1) $|C_n(u) - C_n(v)| \le n^\alpha \rho |u - v|$ for all $u, v \in \mathcal{I}_t$ and some constant $\rho > 0$ depending only on the kernel $k$, the point $t$ and the constant $c$. Also note that both classes $\mathcal{G}_{1,n}$ and $\mathcal{G}_{2,n}$ have a constant $\rho_i$ times $1_{\mathcal{I}_t}$ as envelope function, where the constant $\rho_i$ only depend on $k$, $t$, $c$ and $M_i$, $i = 1, 2$. For $\kappa_{1,n} = n^{3\alpha-5/6}\log n$ and $\kappa_{2,n} = n^{4\alpha-5/6}\log n$, we now have that

$$\begin{aligned} &\mathrm{E}|R_n|1_{\mathcal{E}_{1,n,M_1}} \\ &\quad \le \mathrm{E} \sup_{A\in\mathcal{A}, B\in\mathcal{B}_{M_1}} \left| n^{2\alpha-1/3}\log n \int \psi(A(u))B(u)\, d(\mathbb{G}_n - G)(u) \right| 1_{\mathcal{E}_{1,n,M_1}} \\ &\quad \le \kappa_{1,n} \mathrm{E} \sup_{\xi\in\mathcal{G}_{1,n}} \left| \int \xi(u)\, d\sqrt{n}(\mathbb{G}_n - G)(u) \right| \end{aligned}$$

and

$$\begin{aligned} &\mathrm{E}|S_n|1_{\mathcal{E}_{2,n,M_2}} \\ &\quad \le \mathrm{E} \sup_{B\in\mathcal{B}_{M_2}} \left| n^{2\alpha-1/2} \int h(u,\delta) \right. \end{aligned}$$

$$\times \{\psi(n^{-1/3}B(u) \times \log n + u) - \psi(u)\}\, d\sqrt{n}(\mathbb{P}_n - P_0)(u,\delta)\Bigg| \times 1_{\mathcal{E}_{2,n,M_2}}$$

$$\le \mathrm{E} \sup_{\zeta \in \mathcal{G}_{2,n}} \kappa_{2,n}\left|\int \zeta(u,\delta)\, d\sqrt{n}(\mathbb{P}_n - P_0)(u,\delta)\right|.$$

To bound these expectations, we use Theorem A.6. Using the entropy results for $\mathcal{A}$ and $\mathcal{B}_M$ together with smoothness properties, we bound the entropies of the classes $\mathcal{G}_{1,n}$ and $\mathcal{G}_{2,n}$. Therefore, we fix an arbitrary probability measure $Q$ and $\delta > 0$.

We start with the entropy of $\mathcal{G}_{1,n}$. Select a minimal $n^{-\alpha}\delta/(2\rho M_1)$-net $A_1, \ldots, A_{N_A}$ in $\mathcal{A}$ and a minimal $\delta/(2\|C_n\|_\infty)$-net $B_1, B_2, \ldots, B_{N_B}$ in $\mathcal{B}_{M_1}$ and construct the subset of $\mathcal{G}_{1,n}$ consisting of the functions $\xi_{A_i,B_j,n}$ corresponding to these nets. The number of functions in this net is then given by

$$\begin{aligned} N_A N_B &= \exp(H(n^{-\alpha}\delta/(2\rho M_1), \mathcal{A}, L_2(Q)) + H(\delta/(2\|C_n\|_\infty), \mathcal{B}_{M_1}, L_2(Q))) \\ &\le \exp(Cn^\alpha/\delta), \end{aligned}$$

where $C > 0$ is a constant. This set is a $\delta$-net in $\mathcal{G}_{1,n}$. Indeed, choose a $\xi = \xi_{A,B,n} \in \mathcal{G}_{1,n}$ and denote the closest function to $A$ in the $\mathcal{A}$-net by $A_i$ and similarly the function in the $\mathcal{B}_{M_1}$-net closest to $B$ by $B_j$. Then

$$\begin{aligned} &\|\xi_{A,B,n} - \xi_{A_i,B_j,n}\|_{Q,2} \\ &\quad \le \|C_n\|_\infty \|B(\cdot) - B_j(\cdot)\|_{Q,2} + M_1\|C_n(A_i(\cdot)) - C_n(A(\cdot))\|_{Q,2} \\ &\quad \le \delta/2 + M_1\rho n^\alpha \|A_i - A\|_{Q,2} \le \delta. \end{aligned}$$

This implies that

$$H(\delta, \mathcal{G}_{1,n}, L_2(Q)) \lesssim n^\alpha/\delta$$

and

$$J(\delta, \mathcal{G}_{1,n}) \le \int_0^\delta \sqrt{1 + H(\varepsilon, \mathcal{G}_{1,n}, L_2(Q))}\, d\varepsilon \lesssim n^{\alpha/2}\sqrt{\delta}.$$

To bound the entropy of $\mathcal{G}_{2,n}$, we select a minimal $(\delta/\rho)$-net $B_1, B_2, \ldots, B_N$ in $\mathcal{B}_{M_2}$ and construct the subset of $\mathcal{G}_{2,n}$ consisting of the functions $\zeta_{B_i,n}$ corresponding to this net. The number of functions in this net is then given by

$$N = \exp(H(\delta/\rho, \mathcal{B}_{M_2}, L_2(Q))) \le \exp(C/\delta),$$

where $C>0$ is a constant. This set is a $\delta$-net in $\mathcal{G}_{2,n}$. Indeed, choose a $\zeta=\zeta_{B,n}\in\mathcal{G}_{2,n}$ and denote the closest function to $B$ in the $\mathcal{B}_{M_2}$-net by $B_i$, then

$$\begin{aligned}&\|\zeta_{B,n}-\zeta_{B_i,n}\|_{L_2(Q)}\\&\quad\le n^{1/3-\alpha}(\log n)^{-1}\|h\|_\infty\\&\qquad\times\|C_n(n^{-1/3}B(\cdot)\log n+\cdot)-C_n(n^{-1/3}B_i(\cdot)\log n+\cdot)\|_{L_2(Q)}\\&\quad\le n^{1/3-\alpha}(\log n)^{-1}n^\alpha\rho n^{-1/3}\log n\|B_i-B\|_{L_2(Q)}\le\delta.\end{aligned}$$

This implies that

$$H(\delta,\mathcal{G}_{2,n},L_2(Q))\lesssim 1/\delta\quad\text{and}\quad J(\delta,\mathcal{G}_{2,n})\lesssim\sqrt{\delta}.$$

We now obtain via Theorem A.6 that

$$\begin{aligned}\mathrm{E}|R_n|1_{\mathcal{E}_{1,n,M_1}}&\le\kappa_{1,n}E\sup_{\xi\in\mathcal{G}_{1,n}}\left|\int\xi(u)\,d\sqrt{n}(\mathbb{G}_n-G)(u)\right|\\&\lesssim\kappa_{1,n}J(1,\mathcal{G}_{1,n})\lesssim n^{7\alpha/2-5/6}\log n,\\\mathrm{E}|S_n|1_{\mathcal{E}_{2,n,M_2}}&\le\kappa_{2,n}E\sup_{\zeta\in\mathcal{G}_{2,n}}\left|\int\zeta(u,\delta)\,d\sqrt{n}(\mathbb{P}_n-P_0)(u,\delta)\right|\\&\lesssim\kappa_{2,n}J(1,\mathcal{G}_{2,n})\lesssim n^{4\alpha-5/6}\log n.\end{aligned}$$

Hence, we can take $\nu_{i,n}=\varepsilon n^{-\beta_i}(\log n)^2$ for $\beta_1=5/6-7\alpha/2$, $\beta_2=5/6-4\alpha$ and any $\varepsilon>0$ to conclude that

$$P\left(\frac{n^{\beta_1}}{(\log n)^2}|R_n|>\varepsilon\right)\le\frac{n^{\beta_1}}{\varepsilon(\log n)^2}E|R_n|1_{\mathcal{E}_{1,n,M_1}}+\eta/2\lesssim\frac{1}{\varepsilon\log n}+\eta/2<\eta,$$
$$P\left(\frac{n^{\beta_2}}{(\log n)^2}|S_n|>\varepsilon\right)\le\frac{n^{\beta_2}}{\varepsilon(\log n)^2}E|S_n|1_{\mathcal{E}_{2,n,M_2}}+\eta/2\lesssim\frac{1}{\varepsilon\log n}+\eta/2<\eta$$

for $n$ sufficiently large. □

With this lemma, we now can prove Theorem 4.2.

Proof of Theorem 4.2. Using the piecewise contant version $\bar\psi_{h,t}$ of $\psi_{h,t}$, we can write

$$\int\psi_{h,t}(u)(\delta-\hat F_n(u))\,dP_0(u,\delta)=\int\bar\psi_{h,t}(u)(\delta-\hat F_n(u))\,dP_0(u,\delta)+R_n,$$

where for $h=cn^{-\alpha}$ and $n$ sufficiently large

$$|R_n|\le c_1h^{-2}\int_{u\in[t-h,t+h]}|F_0(u)-\hat F_n(u)|^2\,dG(u)=\mathcal{O}_p(n^{\alpha-2/3})=\mathcal{O}_p(n^{-2\alpha})$$

by (2.3) and Lemma A.4. So we find

$$n^{2\alpha}\int \psi_{h,t}(u)(\delta-\hat F_n(u))\,dP_0(u,\delta)$$
$$=n^{2\alpha}\int \bar\psi_{h,t}(u)(\delta-F_0(u))\,d(P_0-\mathbb{P}_n)(u,\delta)+\mathcal{O}_p(1)$$

using that $n^{2\alpha}R_n=\mathcal{O}_p(1)$, Property A.5 and (A.18). By (A.19), we get

$$n^{2\alpha}\int \bar\psi_{h,t}(u)(\delta-F_0(u))\,d(P_0-\mathbb{P}_n)(u,\delta)$$
$$=n^{2\alpha}\int \psi_{h,t}(u)(\delta-F_0(u))\,d(P_0-\mathbb{P}_n)(u,\delta)+\mathcal{O}_p(1).$$

Applying the central limit theorem with $\alpha=1/5$, gives

$$n^{2/5}\int \psi_{h,t}(u)(\delta-F_0(u))\,d(\mathbb{P}_n-P_0)(u,\delta)\rightsquigarrow \mathcal{N}(0,\sigma^2_{F,\mathrm{SM}})$$

for $\sigma^2_{F,\mathrm{SM}}$ as in (4.4). Note that now

$$\begin{aligned}
&n^{2/5}(\hat F_n^{\mathrm{SM}}(t)-F_0(t))\\
&\quad=n^{2/5}\int \psi_{h,t}(u)(\delta-F_0(u))\,d(\mathbb{P}_n-P_0)(u,\delta)\\
&\qquad+n^{2/5}\left(\int K_h(t-u)\,dF_0(u)-F_0(t)\right)\rightsquigarrow \mathcal{N}(\mu_{F,\mathrm{SM}},\sigma^2_{F,\mathrm{SM}}).
\end{aligned}$$

To find our optimal bandwidth $h_{n,\mathrm{opt}}$, we minimize the aMSE with respect to $c$

$$\mathrm{aMSE}(\hat F_n^{\mathrm{SM}},c)=\frac14 c^4 m_2^2(k) f_0'(t)^2+c^{-1}\frac{F_0(t)(1-F_0(t))}{g(t)}\int k(u)^2\,du,$$

which is standard a minimization in $c$, yielding (4.5). □

**Acknowledgments.** We thank two anonymous referees and the Associate Editor for their valuable remarks.

## REFERENCES

BICKEL, P. J. and FAN, J. (1996). Some problems on the estimation of unimodal densities. *Statist. Sinica* **6** 23–45. MR1379047

DONOHO, D. L. and JOHNSTONE, I. M. (1995). Adapting to unknown smoothness via wavelet shrinkage. *J. Amer. Statist. Assoc.* **90** 1200–1224. MR1379464

DÜMBGEN, L., FREITAG-WOLF, S. and JONGBLOED, G. (2006). Estimating a unimodal distribution function from interval censored data. *J. Amer. Statist. Assoc.* **101** 1094–1106. MR2324149

Dümbgen, L. and Rufibach, K. (2009). Maximum likelihood estimation of a log-concave density and its distribution function: Basic properties and uniform consistency. *Bernoulli* **15** 40–68.

Efron, B. (1979). Bootstrap methods: Another look at the jackknife. *Ann. Statist.* **7** 1–26. MR0515681

Eggermont, P. P. B. and LaRiccia, V. N. (2001). *Maximum Penalized Likelihood Estimation*. Springer, New York. MR1837879

González-Manteiga, W., Cao, R. and Marron, J. S. (1996). Bootstrap selection of the smoothing parameter in nonparametric hazard rate estimation. *J. Amer. Statist. Assoc.* **91** 1130–1140. MR1424613

Groeneboom, P., Jongbloed, G. and Wellner, J. A. (2002). A canonical process for estimation of convex functions: The invelope of integrated brownian motion $+t^4$. *Ann. Statist.* **29** 1620–1652. MR1891741

Groeneboom, P. and Wellner, J. A. (1992). *Information Bounds and Nonparametric Maximum Likelihood Estimation*. Cambridge Univ. Press, New York. MR1180321

Hazelton, M. L. (1996). Bandwidth selection for local density estimators. *Scand. J. Statist.* **23** 221–232. MR1394655

Huang, Y. and Zhang, C. H. (1994). Estimating a monotone density from censored observations. *Ann. Statist.* **22** 1256–1274. MR1311975

Jones, M. C. (1993). Simple boundary correction for kernel density estimation. *Stat. Comput.* **3** 135–146.

Kaplan, E. L. and Meier, P. (1958). Nonparametric estimation from incomplete data. *J. Amer. Statist. Assoc.* **53** 457–481. MR0093867

Keiding, N. (1991). Age-specific incidence and prevalence: A statistical perspective. *J. Roy. Statist. Soc. Ser. A* **154** 371–412. MR1144166

Mammen, E. (1991). Estimating a monotone regression function. *Ann. Statist.* **19** 724–740. MR1105841

Marron, J. S. and Padgett, W. J. (1987). Asymptotically optimal bandwidth selection for kernel density estimators from randomly right-censored samples. *Ann. Statist.* **15** 1520–1535. MR0913571

Patil, P. N., Wells, M. T. and Marron, J. S. (1994). Some heuristics of kernel based estimators of ratio functions. *J. Nonparametr. Stat.* **4** 203–209. MR1290930

Schuster, E. F. (1985). Incorporating support constraints into nonparametric estimators of densities. *Comm. Statist. Theory Methods* **14** 1123–1136. MR0797636

Sheather, S. J. (1983). A data-based algorithm for choosing the window width when estimating the density at a point. *Comput. Statist. Data Anal.* **1** 229–238.

Silverman, B. W. (1978). Weak and strong uniform consistency of the kernel estimate of a density and its derivative. *Ann. Statist.* **6** 177–184. MR0471166

Silverman, B. W. (1986). *Density Estimation for Statistics and Data Analysis*. Chapman and Hall, London. MR0848134

van de Geer, S. A. (2000). *Applications of Empirical Process Theory*. Cambridge Univ. Press, New York. MR1739079

van der Vaart, A. W. and van der Laan, M. J. (2003). Smooth estimation of a monotone density. *Statistics* **37** 189–203. MR1986176

van der Vaart, A. W. and Wellner, J. A. (1996). *Weak Convergence and Empirical Processes*. Springer, New York. MR1385671

P. Groeneboom
G. Jongbloed
B. I. Witte
Delft University of Technology
Delft Institute of Applied Mathematics
Mekelweg 4
2628 CD Delft
The Netherlands
E-mail: P.Groeneboom@tudelft.nl
G.Jongbloed@tudelft.nl
B.I.Witte@tudelft.nl